\documentclass[12pt,reqno]{amsart}

\usepackage{amssymb}
\usepackage{dsfont,enumerate,float}
\usepackage{tikz}
\usepackage[margin=2.2cm]{geometry}

\usepackage[T1]{fontenc}
\usepackage[english]{babel}
\usepackage{lmodern,microtype,hyperref}

\allowdisplaybreaks
\hypersetup{hidelinks}

\DeclareMathOperator{\curl}{curl}
\DeclareMathOperator{\Div}{div}
\DeclareMathOperator{\dist}{dist}
\DeclareMathOperator{\Op}{Op}
\DeclareMathOperator{\Rer}{Re}
\DeclareMathOperator{\Vol}{Vol}

\newcommand{\bs}{\backslash}
\newcommand{\C}{\mathbb{C}}
\newcommand{\D}{\mathbb{D}}
\newcommand{\Dc}{\mathcal{D}}
\newcommand{\GL}{\mathrm{GL}}
\newcommand{\Hb}{\mathbb{H}}
\newcommand{\indicator}{\mathds{1}}
\newcommand{\M}{\mathbb{M}}
\newcommand{\N}{\mathbb{N}}
\newcommand{\ol}{\overline}
\newcommand{\Pc}{\mathcal{P}}
\newcommand{\PSL}{\mathrm{PSL}}
\newcommand{\R}{\mathbb{R}}
\newcommand{\Sb}{\mathbb{S}}

\newcommand{\T}{\mathbb{T}}
\newcommand{\ve}{\varepsilon}
\newcommand{\vp}{\varphi}
\newcommand{\Z}{\mathbb{Z}}

\theoremstyle{plain}
\newtheorem{thm}{Theorem}[section]

\newtheorem{prob}{Problem}
\newtheorem{lemma}[thm]{Lemma}
\newtheorem{prop}[thm]{Proposition}

\theoremstyle{definition}
\newtheorem*{defn}{Definition}
\newtheorem*{ex}{Example}
\newtheorem*{rmk}{Remark}

\numberwithin{equation}{section}

\title[Invariant measures with uniform projection]{Invariant measures with uniform projection on surfaces and billiards}
\author{Xiaolong Han}
\email{xiaolong.han@csun.edu}
\address{Department of Mathematics, California State University, Northridge, CA 91330, USA}


\begin{document}

\begin{abstract}
In this paper, we study probability measures on the unit tangent bundle which are invariant under the geodesic flow and whose projection to the base manifold is the normalized Riemannian volume. We show that the Liouville measure is the only such measure on a complete oriented surface of finite area with ergodic geodesic flow. In particular, this resolves a conjecture of Katok--Knieper--Weiss in dimension two. The same holds for ergodic planar billiards with piecewise $C^2$ boundary, such as the Bunimovich stadium, the Sinai billiard, and generic polygons. The proof combines the Fourier analysis in the fibers with holomorphic functions on the transport twistor space.

As a consequence, on these surfaces and on the Bunimovich stadium, a high-energy sequence of normalized Laplacian eigenfunctions tends to be equidistributed in the phase space if and only if it tends to be equidistributed in the physical space. In particular, on every stadium, every failure of equidistribution in the phase space implies scarring in the physical space. Combined with Hassell's theorem, this gives a new proof of scarring in the physical space for almost every stadium.

In higher dimensions, we obtain an entropy bound for the ergodic components of such measures. We also construct counterexamples for hyperbolic toral automorphisms.
\end{abstract}

\maketitle

\setcounter{tocdepth}{1}
\tableofcontents

%%%%%%%%%%%%%%%%%%%%%%%%%%%%%%%%%%%%%%%%%
\section{Introduction}\label{sec:intro}
Let $\M$ be a connected complete Riemannian manifold of finite volume. Denote by $S\M$ its unit tangent bundle and by $a_t$ the geodesic flow on $S\M$. Write $\Pi:S\M\to\M$ for the natural projection. Throughout the paper, $d\Vol_\M$ and $d\Vol_L$ denote the normalized Riemannian volume on $\M$ and the normalized Liouville measure on $S\M$, respectively. Thus,
\[\Pi_*d\Vol_L=d\Vol_\M.\]
We say that a probability measure $\mu$ on $S\M$ has uniform projection if $\Pi_*\mu=d\Vol_\M$. We consider the following problem.
\begin{prob}[Invariant measures with uniform projection]\label{prob:uniform}
Suppose that $a_t$ is ergodic with respect to $d\Vol_L$. Let $\mu$ be a probability measure on $S\M$ such that
\begin{equation}\label{eq:IP}
\left(a_t\right)_*\mu=\mu\text{ for }t\in\R\quad\text{and}\quad\Pi_*\mu=d\Vol_\M.
\end{equation}
Determine whether $\mu=d\Vol_L$.
\end{prob}

\begin{ex}
The ergodicity hypothesis in Problem \ref{prob:uniform} cannot be dropped. For example, let $\T^n=\R^n/\Z^n$ with $n\ge2$ be the flat torus equipped with the Lebesgue measure $dx$. Fix a unit vector $v_0\in\R^n$. Then the measure $dx\otimes\delta_{v_0}$ on $S\T^n=\T^n\times\Sb^{n-1}$ is invariant under the geodesic flow. Its projection to $\T^n$ is $dx$, but it is not the Liouville measure. Here, the geodesic flow is not ergodic with respect to the Liouville measure. In Section \ref{sec:surfaceproof}, we use this example to show where the ergodicity enters the proof of Theorem \ref{thm:surface}.

The analogous statement can fail for a discrete-time dynamical system, even when the reference measure is ergodic. For example, let $M:\T^2\to\T^2$ be a hyperbolic automorphism of the two-dimensional torus $\T^2=\R^2/\Z^2$, regarded as the phase space of the one-dimensional torus $\T=\R/\Z$. These automorphisms are the classical maps underlying quantum cat maps; see Faure--Nonnenmacher--De Bi\`evre \cite{FNDB}. Every such $M$ is ergodic with respect to the Lebesgue measure on $\T^2$. In Appendix \ref{sec:toral}, we construct, for every such $M$, an $M$-invariant probability measure whose projection to the position is the Lebesgue measure on $\T$, but which is singular with respect to the Lebesgue measure on $\T^2$.
\end{ex}

This problem is motivated by the study of Laplacian eigenfunctions in Quantum Chaos. In particular, the quantum unique ergodicity (QUE) conjecture of Rudnick--Sarnak \cite{RS} asserts that on compact negatively curved manifolds, every high-energy sequence of normalized Laplacian eigenfunctions tends to be equidistributed in the phase space. In this setting, the geodesic flow is ergodic by the theorems of Hopf \cite{Ho1, Ho2} and Anosov \cite{Ano}. The semiclassical measures describe the distribution of eigenfunctions in the phase space $S\M$, while their projections describe the distribution in the physical space $\M$. Equidistribution in the phase space implies equidistribution in the physical space. A positive answer to Problem \ref{prob:uniform} implies the converse for eigenfunction sequences; see Section \ref{sec:QUE}.

In fact, Katok--Knieper--Weiss \cite[Conjecture 3]{KKW} conjectured that on a compact negatively curved manifold, the Liouville measure is the only invariant probability measure with uniform projection. Their conjecture arose from their study of critical points of topological entropy under perturbations of the Riemannian metric. In particular, they proved that on compact surfaces, these critical points among negatively curved metrics of fixed area are exactly the metrics of constant negative curvature.

Several subsequent works studied related questions. Flaminio \cite{F} proved that even invariant distributions on the unit tangent bundle of a positively curved sphere are determined by their projections. Ledrappier--Lindenstrauss \cite[Theorem 1.1]{LeLi} considered an arbitrary invariant probability measure $\mu$ on the unit tangent bundle $S\M$ of a compact surface $\M$. With $\dim\mu$ defined as the essential infimum of the pointwise dimensions, they proved that, among other things, if $\dim\mu\le2$, then $\dim(\Pi_*\mu)=\dim\mu$, whereas if $\dim\mu>2$, then $\Pi_*\mu\ll d\Vol_\M$. These results concern the dimension and the absolute continuity of the projection. Their theorem does not require negative curvature, and it does not determine the invariant measures whose projection equals $d\Vol_\M$.

Hovila--J\"arvenp\"a\"a--J\"arvenp\"a\"a--Ledrappier \cite{HJJL} constructed ergodic invariant measures on compact negatively curved surfaces whose projections have Hausdorff dimension $2$ but are singular with respect to the Riemannian area. Thus, full dimension of the projection does not imply uniform projection.

In this paper, we resolve Conjecture 3 of \cite{KKW} in dimension $2$. More generally, we establish the conclusion for complete oriented surfaces of finite area whose geodesic flow is ergodic with respect to the Liouville measure. Our first main theorem is the following.
\begin{thm}[Invariant measures with uniform projection on surfaces]\label{thm:surface}
Let $\M$ be a connected complete oriented Riemannian surface of finite area. Suppose that the geodesic flow is ergodic with respect to $d\Vol_L$. If a probability measure $\mu$ on $S\M$ satisfies \eqref{eq:IP}, then $\mu=d\Vol_L$.
\end{thm}

It is well known that the geodesic flow is ergodic on compact negatively curved surfaces and on finite-area hyperbolic surfaces; see Hedlund \cite{Hed} and Hopf \cite{Ho1, Ho2}. Therefore, Theorem \ref{thm:surface} applies to these surfaces. For example, on a compact negatively curved surface, the measure of maximal entropy has uniform projection if and only if the curvature is constant. Indeed, by Katok \cite{K}, the measure of maximal entropy is the Liouville measure if and only if the curvature is constant. In a recent preprint, Rodriguez Hertz--Rodriguez Hertz--Wu \cite{RRW} proved that the geodesic flow is ergodic on closed oriented surfaces of genus at least two with non-positive curvature. So Theorem \ref{thm:surface} also applies to these surfaces. The conclusion also extends to compact non-orientable negatively curved surfaces and to finite-area hyperbolic orbifolds, including the modular surface, by passage to finite covers (Proposition \ref{prop:lattice}).

The proof uses the Fourier analysis of invariant distributions in the fibers of $S\M$, which goes back to Guillemin--Kazhdan \cite{GK} and was developed by Paternain--Salo--Uhlmann \cite{PSU1, PSU2}, and the transport twistor space introduced by Bohr--Paternain \cite{BP} and developed by Bohr--Lefeuvre--Paternain \cite{BLP}; see also Bohr--Monard--Paternain \cite{BMP}. The transport twistor space is the interior $Z^\circ$ of the closed unit ball bundle
\begin{equation}\label{eq:Z}
Z=\left\{(x,w)\in T\M:|w|_x\le1\right\},
\end{equation}
equipped with a complex structure relating fiberwise holomorphic invariant distributions on $S\M$ to holomorphic functions on $Z^\circ$ with controlled growth at the boundary; see Section \ref{sec:twistor}. A similar approach based on the Fourier analysis in the fibers was recently used in \cite{RRW} to solve the long-standing problem of the ergodicity of the geodesic flow on surfaces of non-positive curvature.

\begin{rmk}[Outline of the proof]
We first average $\mu$ with its image under the velocity reversal. The average is still invariant under $a_t$ and has uniform projection, and its odd Fourier modes in the fibers vanish. The invariance gives a recurrence relation between the Fourier modes. Since the zeroth Fourier mode is equal to $1$, the recurrence relation also holds for the distribution obtained by retaining the zeroth Fourier mode and doubling the positive Fourier modes. Hence, this distribution extends to a holomorphic function $H$ on $Z^\circ$. The restriction of $H$ to each fiber of $Z^\circ$ is the Carath\'eodory transform of the conditional measure on the corresponding unit circle. Since the conditional measures are probability measures, we have $\Rer H>0$ and $H=1$ on the zero section. So the Cayley transform of $H$ is a bounded holomorphic function on $Z^\circ$. Its boundary value is invariant under the geodesic flow, and the ergodicity forces it to be constant. Therefore, the average of $\mu$ is the Liouville measure. Finally, the Liouville measure is an extreme point of the invariant probability measures by the ergodicity, so $\mu$ is the Liouville measure.
\end{rmk}

The method also applies to billiards, but with necessary modifications due to the boundary. Let $\Omega\subset\R^2$ be a bounded connected open set. Suppose that $\partial\Omega$ is the union of finitely many regular embedded $C^2$ curves, which meet only at their endpoints. Denote by $\Sigma$ the finite set of these endpoints. Suppose that $\Omega$ lies on one side of $\partial\Omega$ near each point of $\partial\Omega\setminus\Sigma$. Identify the unit tangent bundle of $\ol\Omega$ with $\ol\Omega\times\Sb$, where the point $(x,\theta)$ represents the direction $(\cos\theta,\sin\theta)$ at $x$.

Still write $\Pi$ for the projection $\ol\Omega\times\Sb\to\ol\Omega$. We say that a point of $\Sigma$ is a corner if $\partial\Omega$ is not a $C^1$ curve near it. For example, the corners of a polygon are its vertices. Denote by $a_t$ the billiard flow, that is, a point moves along a straight line at unit speed, and its direction is reflected across the tangent line when it hits $\partial\Omega$ at a point which is not a corner. The billiard flow is defined at the points whose trajectories do not hit the corners and have finitely many reflections in every bounded time interval. These points have full measure with respect to $d\Vol_L$; see Chernov--Markarian \cite[Chapter 2]{CM}.

A probability measure $\mu$ on $\ol\Omega\times\Sb$ is called invariant if $a_t$ is defined $\mu$-almost everywhere and $(a_t)_*\mu=\mu$ for $t\in\R$. Let $dx$ and $d\theta$ be the Lebesgue measures on $\R^2$ and on $\Sb$, respectively. Write $|\Omega|$ for the Lebesgue measure of $\Omega$. Then
\[d\Vol_\Omega=\frac{dx}{|\Omega|}\quad\text{and}\quad d\Vol_L=\frac{dxd\theta}{2\pi|\Omega|}.\]
Our second main theorem states that Problem \ref{prob:uniform} has a positive answer for ergodic billiards.
\begin{thm}[Invariant measures with uniform projection on billiards]\label{thm:billiard}
Let $\Omega$ be as above. Suppose that the billiard flow is ergodic with respect to $d\Vol_L$. If an invariant probability measure $\mu$ on $\ol\Omega\times\Sb$ satisfies $\Pi_*\mu=d\Vol_\Omega$, then $\mu=d\Vol_L$.
\end{thm}

\begin{rmk}[Outline of the proof]
In the interior of $\Omega$, the proof is the same as that of Theorem \ref{thm:surface}. It gives a bounded holomorphic function $F$ on the transport twistor space. For the Euclidean metric, the transport twistor space has the explicit holomorphic coordinates $(w,z-w^2\bar z)$, where $z=x_1+ix_2$ and $w\in\D$ is the complex coordinate in the fibers.

The difference is in the use of the ergodicity. On a surface without boundary, the limit $f$ of $F$ at the boundary of the fibers is invariant under the geodesic flow. On a billiard table, $f$ is only invariant along the straight segments of the trajectories. The reflection at $\partial\Omega$ reverses the orientation of the fibers. So it maps the fiberwise holomorphic functions to fiberwise anti-holomorphic functions, and $f$ is in general not invariant under the billiard flow.

The boundary is treated in two steps. First, $F$ has boundary values $F^\partial$ on the fibers over $\partial\Omega$. Write $\tau_x$ for the reflection across the tangent line at $x$. The invariance of $\mu$ under the reflections gives
\[F^\partial\left(x,\tau_x(w)\right)=\ol{F^\partial(x,w)}.\]
Hence, the real parts of these boundary values are invariant under reflection. Second, this symmetry implies that $\Rer f$ is invariant under the billiard flow. The Fourier modes of $f$ of positive weight are determined by those of $\Rer f$. So the ergodicity still forces $F$ to be constant.
\end{rmk}

Theorem \ref{thm:billiard} applies to the following ergodic billiards; see \cite{CM} for a complete treatment of billiard dynamics.

\begin{figure}[!ht]
\centering
\begin{tikzpicture}
\draw[thick] (-2.5,-2.5) rectangle (2.5,2.5);
\draw[thick] (0,0) circle (1.2);
\draw[|-|] (-2.5,-2.9) -- (2.5,-2.9) node[midway,below] {$[0,1]$};
\draw[|-|] (2.9,-2.5) -- (2.9,2.5) node[midway,right] {$[0,1]$};
\draw[->,thick] (-2.193,-0.732) -- (-1.186,-0.182);
\draw[->,thick] (-1.186,-0.182) -- (-2.5,0.078);
\draw[->,thick] (-2.5,0.078) -- (-1.15,0.344);
\draw[->,thick] (-1.15,0.344) -- (-1.791,0.976);
\node at (1.7,1.8) {$\Omega$};
\end{tikzpicture}
\caption{The Sinai billiard $\Omega$, a square with a disk removed. A trajectory of the billiard flow reflects at the disk, at a side, and at the disk again.}\label{fig:sinai}
\end{figure}
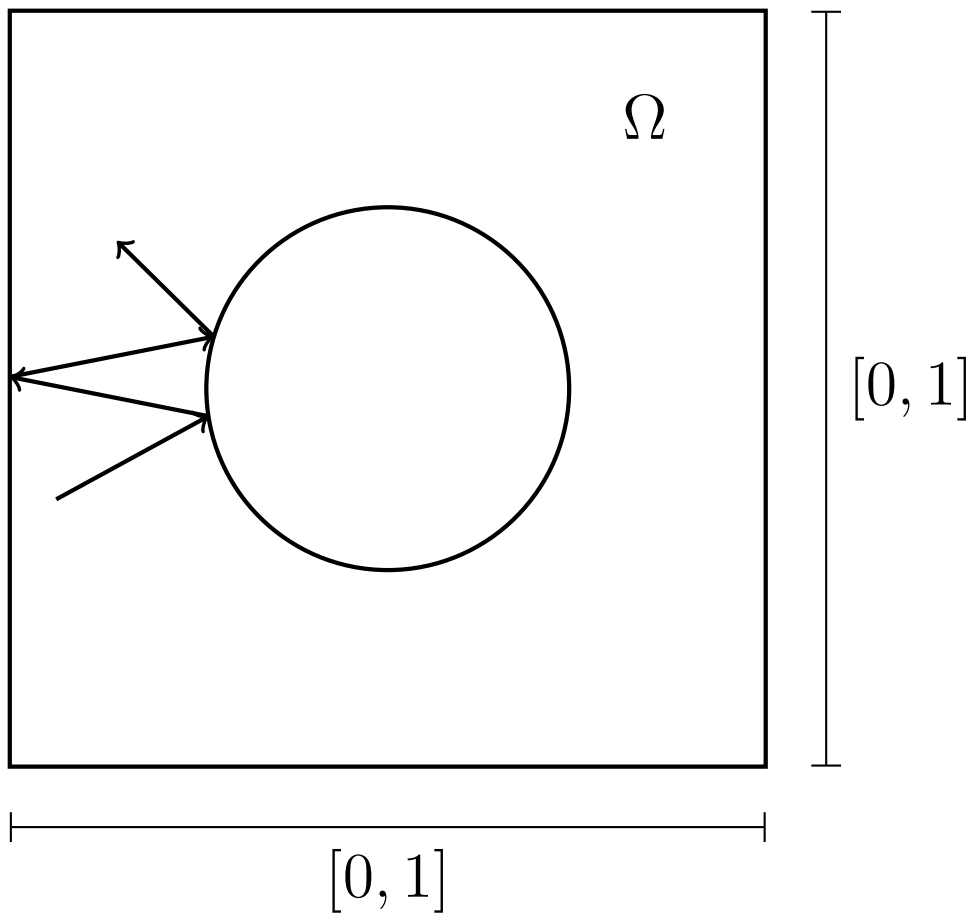

$\bullet$ The Sinai billiard. The table $\Omega$ is a square with a closed disk removed from its center; see Figure \ref{fig:sinai}. Sinai \cite{Sin} proved that the billiard flow is ergodic.

$\bullet$ The Bunimovich stadia. They are defined in Section \ref{sec:scar}; see Figure \ref{fig:stadium}. Bunimovich \cite{B} proved that the billiard flow is ergodic.

$\bullet$ The generic polygons. Here, $\Sigma$ is the set of the vertices. Kerckhoff--Masur--Smillie \cite{KMS} proved that the billiard flow is ergodic for a dense $G_\delta$ set of polygons with a given number of vertices.

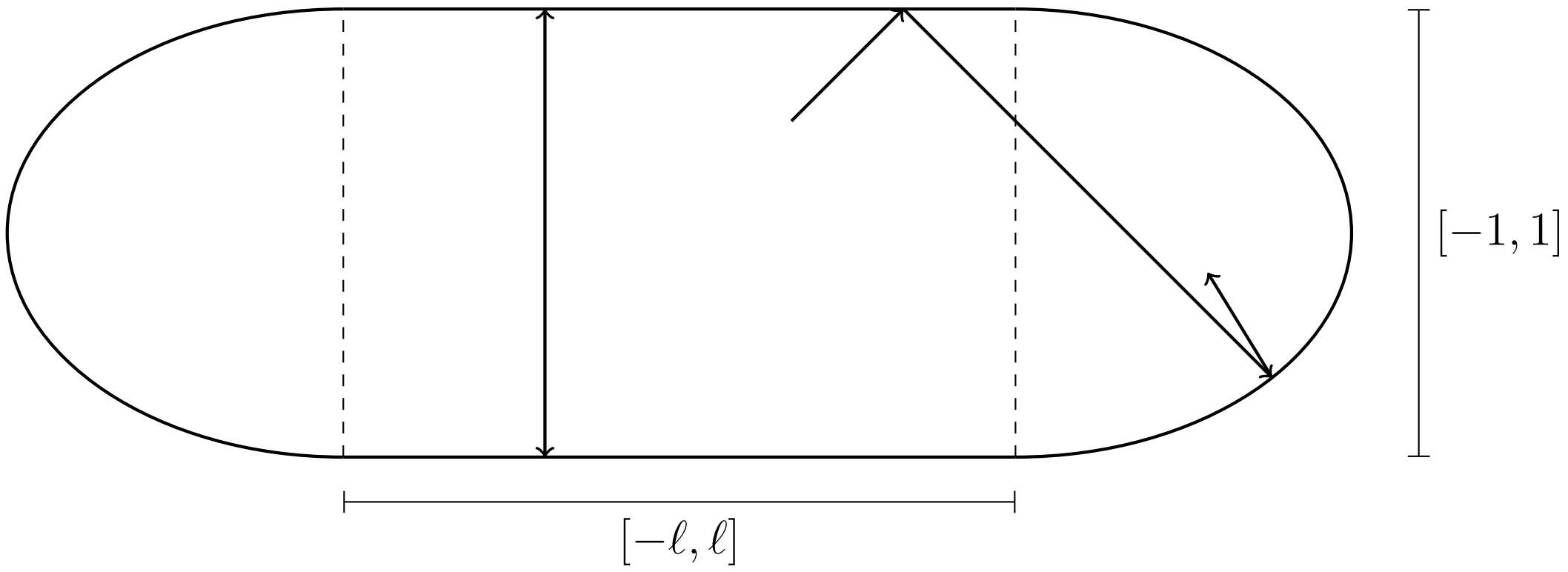
\begin{figure}[!ht]
\centering
\begin{tikzpicture}
\draw[thick] (-3,2) -- (3,2);
\draw[thick] (-3,-2) -- (3,-2);
\draw[thick] (3,-2) arc[start angle=270,end angle=450,x radius=3,y radius=2];
\draw[thick] (-3,2) arc[start angle=90,end angle=270,x radius=3,y radius=2];
\draw[dashed] (-3,-2) -- (-3,2);
\draw[dashed] (3,-2) -- (3,2);
\draw[|-|] (-3,-2.4) -- (3,-2.4) node[midway,below] {$[-\ell,\ell]$};
\draw[|-|] (6.6,-2) -- (6.6,2) node[midway,right] {$[-1,1]$};
\draw[->,thick] (1,1) -- (2,2);
\draw[->,thick] (2,2) -- (5.291,-1.291);
\draw[->,thick] (5.291,-1.291) -- (4.716,-0.354);
\draw[<->,thick] (-1.2,-2) -- (-1.2,2);
\end{tikzpicture}
\caption{An elliptical stadium $\Omega$. A trajectory of the billiard flow reflects at the boundary, and a bouncing ball trajectory is vertical in the rectangle.}\label{fig:elliptic}
\end{figure}

$\bullet$ The elliptical stadia. Let $\ell>0$ and $b>0$. The table $\Omega$ is the interior of the union of the closed rectangle $[-\ell,\ell]\times[-1,1]$ and two closed half-ellipses attached to its vertical sides. The half-ellipses have vertical semi-axis $1$ and horizontal semi-axis $b$; see Figure \ref{fig:elliptic}. For $b=1$, this is the Bunimovich stadium. Markarian--Oliffson Kamphorst--Pinto de Carvalho \cite[Theorem 2]{MKC} proved that the billiard flow is ergodic if
\[1<b<\sqrt{4-2\sqrt2}\quad\text{and}\quad\ell>2b^2\sqrt{b^2-1}.\]
These tables were recently considered by Li--Zworski \cite{LiZw}; see Section \ref{sec:scar}.

\begin{rmk}
The ergodicity hypothesis cannot be dropped in Theorem \ref{thm:billiard}. Let $\Omega$ be a rational polygon, that is, a polygon whose angles are rational multiples of $\pi$. Denote by $W$ the group generated by the reflections of $\Sb$ across the directions of the sides. Since the angles are rational multiples of $\pi$, the group $W$ is finite. For $\theta\in\Sb$, the measure
\[d\Vol_\Omega\otimes\frac1{|W\theta|}\sum_{\alpha\in W\theta}\delta_\alpha\]
is invariant under the billiard flow and has uniform projection, but it is not the Liouville measure. On the quantum side, Marklof--Rudnick \cite{MR} proved that every orthonormal basis of Dirichlet eigenfunctions on a rational polygon has a density-one subsequence which tends to be equidistributed in the physical space. Along this density-one subsequence, any scarring can occur only in the momentum variables.
\end{rmk}

In dimensions $n\ge3$, Problem \ref{prob:uniform} remains open. Our third main theorem gives a restriction on any possible counterexample. For an invariant probability measure $\nu$, denote by $h_\nu(a_1)$ the measure-theoretic entropy of the time-one map $a_1$ with respect to $\nu$.
\begin{thm}[Entropy of the ergodic components]\label{thm:entropy}
Let $\M$ be a compact connected Riemannian manifold of dimension $n\ge2$. Suppose that its sectional curvature is at most $-K^2$ for some $K>0$. If a probability measure $\mu$ on $S\M$ satisfies \eqref{eq:IP}, then
\[h_{\mu_z}\left(a_1\right)\ge\frac{(n-1)K}{2}\quad\text{for }\mu\text{-almost every }z\in S\M,\]
where $\mu_z$ denotes the component at $z$ in the ergodic decomposition of $\mu$ with respect to the geodesic flow.
\end{thm}

If $\M$ has constant sectional curvature $-K^2$, then the lower bound in Theorem \ref{thm:entropy} is half the topological entropy $(n-1)K$ of $a_1$. The conclusion is component-wise, although the ergodic components of $\mu$ need not have uniform projection. The proof only uses the bound $\Pi_*\mu\le C\,d\Vol_\M$ for some $C>0$.

\subsection{Equidistribution of eigenfunctions on surfaces}\label{sec:QUE}
Quantum Chaos concerns the behavior of quantum systems whose classical dynamics is chaotic. For the geodesic flow on $\M$, the quantum states are the Laplacian eigenfunctions. Write $\Delta$ for the non-negative Laplacian on $\M$. Let $J\subset\N$ be an infinite set, where $\N=\{0,1,2,\dots\}$. Suppose that $\{\vp_j\}_{j\in J}$ is a sequence of eigenfunctions such that
\[\Delta\vp_j=\lambda_j^2\vp_j,\quad\text{where }\|\vp_j\|_{L^2(\M)}=1\text{ and }\lambda_j\to\infty\text{ as }j\in J,\ j\to\infty.\]
All $L^2(\M)$ norms and inner products are taken with respect to $d\Vol_\M$. Write $h_j=\lambda_j^{-1}$ for the semiclassical parameter. Use a semiclassical quantization $\Op_h$ of the symbols in $C_0^\infty(T^*\M)$, chosen such that $\Op_h(f)$ has a compactly supported Schwartz kernel; see Zworski \cite[Section 14.2]{Zw}. Identify $S^*\M$ with $S\M$ by the metric. For $f\in C_0^\infty(T^*\M)$, the matrix element
\[\left\langle\Op_{h_j}(f)\vp_j,\vp_j\right\rangle_{L^2(\M)}\]
describes the distribution of the eigenfunction $\vp_j$ in the phase space.

It is well known that every sequence of eigenfunctions has a subsequence along which these matrix elements converge for all phase-space observables $f\in C_0^\infty(T^*\M)$. The limit is the integral of $f$ over $S\M$ against a non-negative measure $\nu$ on $S\M$, which is called a semiclassical measure and is invariant under $a_t$; see \cite[Sections 5.1 and 5.2]{Zw}. If $\M$ is non-compact, then $\nu(S\M)$ may be less than $1$ because of the escape of mass. For physical-space observables $f\in C_0^\infty(\M)$, the multiplication by $f$ agrees microlocally near $S^*\M$ with the quantization of a compactly supported symbol equal to $f\circ\Pi$ near $S^*\M$. Therefore, the limits of the measures $|\vp_j|^2\,d\Vol_\M$ along the same subsequence are given by $\Pi_*\nu$.

\begin{defn}[Equidistribution of eigenfunctions]\leavevmode
\begin{itemize}
\item We say that $\{\vp_j\}_{j\in J}$ tends to be \textit{equidistributed in the physical space} if
\[\lim_{j\in J,j\to\infty}\int_\M f|\vp_j|^2\,d\Vol_\M=\int_\M f\,d\Vol_\M\quad\text{for every }f\in C_0^\infty(\M).\]
\item We say that $\{\vp_j\}_{j\in J}$ tends to be \textit{equidistributed in the phase space} if
\[\lim_{j\in J,j\to\infty}\left\langle\Op_{h_j}(f)\vp_j,\vp_j\right\rangle_{L^2(\M)}=\int_{S\M}f\,d\Vol_L\quad\text{for every }f\in C_0^\infty(T^*\M).\]
\end{itemize}
\end{defn}

On compact manifolds with ergodic geodesic flow, the quantum ergodicity theorem of Shnirelman \cite{Shn}, Zelditch \cite{Ze}, and Colin de Verdi\`ere \cite{CdV} asserts that every orthonormal basis of eigenfunctions has a density-one subsequence which tends to be equidistributed in the phase space. The QUE conjecture of Rudnick--Sarnak \cite{RS} asserts that on a compact negatively curved manifold, every sequence of normalized eigenfunctions with eigenvalues tending to infinity tends to be equidistributed in the phase space. It remains open; see Dyatlov \cite{D} for an overview of related progress.

Anantharaman \cite{Ana} proved that semiclassical measures on compact manifolds with Anosov geodesic flow have positive entropy. Anantharaman--Nonnenmacher \cite{AN} established a quantitative entropy lower bound which, in constant sectional curvature $-K^2$, equals $(n-1)K/2$. On compact Anosov surfaces, Rivi\`ere \cite{Ri} proved the lower bound given by half the integral of the positive Lyapunov exponent. In constant negative curvature, these semiclassical entropy bounds have the same numerical value as the bound in Theorem \ref{thm:entropy}, but the hypotheses and conclusions are different. These results concern semiclassical measures of Laplacian eigenfunctions, whereas Theorem \ref{thm:entropy} applies to arbitrary invariant probability measures with uniform projection and gives a lower bound for almost every ergodic component.

Theorem \ref{thm:surface} has the following consequence.
\begin{thm}[Equidistribution in the physical space and in the phase space]\label{thm:QUE}
Let $\M$ be as in Theorem \ref{thm:surface}. Suppose that $\{\vp_j\}_{j\in J}$ is a sequence of normalized eigenfunctions with eigenvalues tending to infinity. Then the following are equivalent.
\begin{enumerate}[(i)]
\item $\{\vp_j\}_{j\in J}$ tends to be equidistributed in the physical space.
\item $\{\vp_j\}_{j\in J}$ tends to be equidistributed in the phase space.
\end{enumerate}
\end{thm}

\begin{proof}
(ii)$\Rightarrow$(i). For $f\in C_0^\infty(\M)$, apply (ii) to a symbol in $C_0^\infty(T^*\M)$ which equals $f\circ\Pi$ near $S^*\M$.

(i)$\Rightarrow$(ii). Let $\nu$ be the semiclassical measure of a subsequence. By (i), $\Pi_*\nu=d\Vol_\M$. In particular, $\nu(S\M)=1$. That is, there is no escape of mass. By Theorem \ref{thm:surface}, $\nu=d\Vol_L$. Since every subsequence has a further subsequence with this limit, (ii) holds.
\end{proof}

Therefore, on these surfaces, QUE follows if every sequence of normalized eigenfunctions with eigenvalues tending to infinity tends to be equidistributed in the physical space.

\subsection{Scarring of eigenfunctions in the stadium}\label{sec:scar}
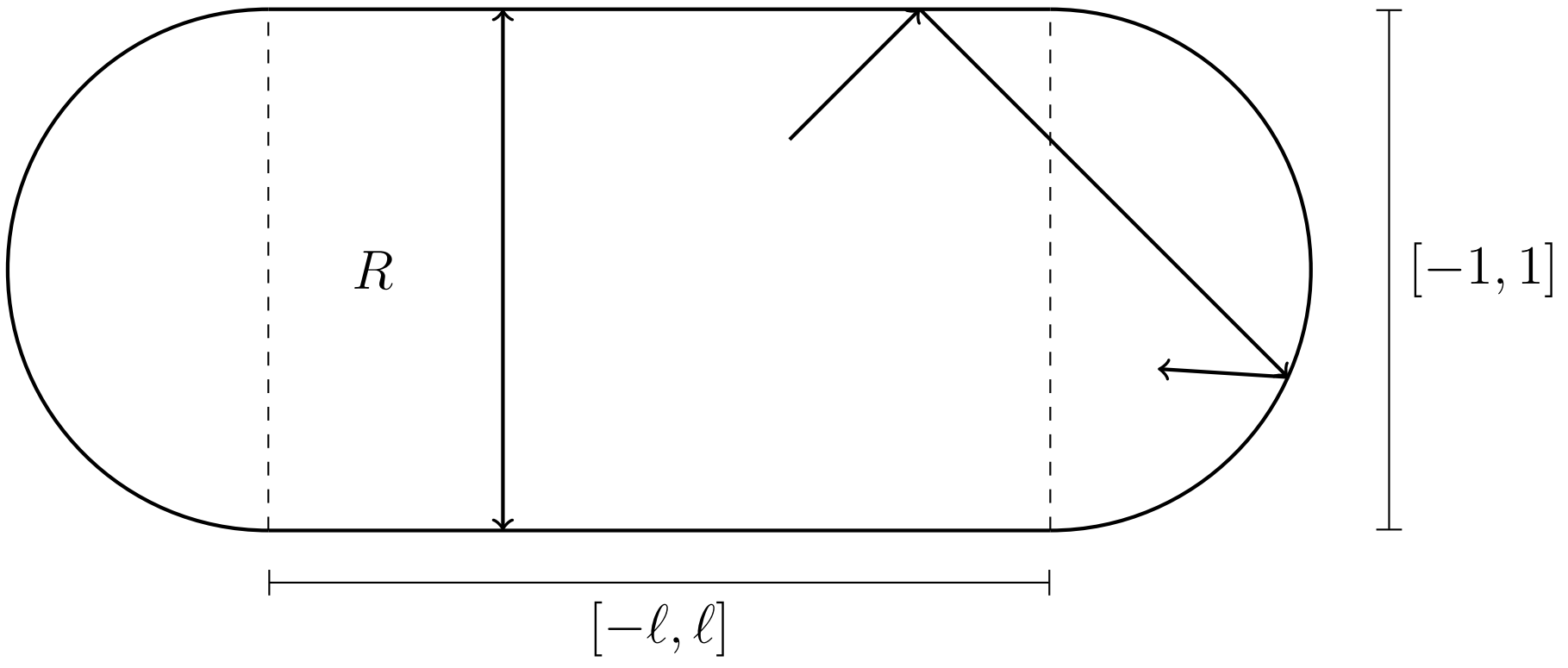
\begin{figure}[!ht]
\centering
\begin{tikzpicture}
\draw[thick] (-3,2) -- (3,2);
\draw[thick] (-3,-2) -- (3,-2);
\draw[thick] (3,-2) arc (270:450:2);
\draw[thick] (-3,2) arc (90:270:2);
\draw[dashed] (-3,-2) -- (-3,2);
\draw[dashed] (3,-2) -- (3,2);
\draw[|-|] (-3,-2.4) -- (3,-2.4) node[midway,below] {$[-\ell,\ell]$};
\draw[|-|] (5.6,-2) -- (5.6,2) node[midway,right] {$[-1,1]$};
\draw[->,thick] (1,1) -- (2,2);
\draw[->,thick] (2,2) -- (4.823,-0.823);
\draw[->,thick] (4.823,-0.823) -- (3.825,-0.76);
\draw[<->,thick] (-1.2,-2) -- (-1.2,2);
\node at (-2.2,0) {$R$};
\end{tikzpicture}
\caption{The Bunimovich stadium $\Omega$ with the rectangle $R$. A trajectory of the billiard flow reflects at the boundary, and a bouncing ball trajectory is vertical in $R$.}\label{fig:stadium}
\end{figure}

In this subsection, we consider the consequences of Theorem \ref{thm:billiard} for the distribution of eigenfunctions on the Bunimovich stadium. Let $\ell>0$. Write
\[R=(-\ell,\ell)\times(-1,1).\]
The Bunimovich stadium $\Omega$ is the interior of the union of the closed rectangle $\ol R$ and the two closed half-disks of radius $1$ attached to its vertical sides; see Figure \ref{fig:stadium}.

Consider the Laplacian on $\Omega$ with the Dirichlet boundary condition. Other boundary conditions, such as the Neumann and Robin boundary conditions, can in principle be treated similarly, but we do not pursue them for simplicity of presentation.

Define the semiclassical measures and the equidistribution of eigenfunctions as in Section \ref{sec:QUE}, with $\M$ replaced by $\Omega$. A semiclassical measure $\nu$ is a measure on $\Omega\times\Sb$ of total mass at most $1$. The total mass may be less than $1$, since part of the mass of the eigenfunctions may concentrate at $\partial\Omega$. The stadium is convex, and its boundary is $C^{1,1}$ and has no corners. So the billiard flow is defined at every point of $\Omega\times\Sb$. By G\'erard--Leichtnam \cite[Theorem 2.3 and Corollary 2.6]{GL}, which apply to convex domains with $C^{1,1}$ boundary, $\nu$ is invariant under the billiard flow. We say that $\Omega$ is QUE if every sequence of normalized eigenfunctions with eigenvalues tending to infinity tends to be equidistributed in the phase space.

Theorem \ref{thm:billiard} has the following consequence.
\begin{thm}[Equidistribution in the stadium]\label{thm:QUEstadium}
Let $\Omega$ be a Bunimovich stadium with the Dirichlet boundary condition. Suppose that $\{\vp_j\}_{j\in J}$ is a sequence of normalized eigenfunctions with eigenvalues tending to infinity. Then the following are equivalent.
\begin{enumerate}[(i)]
\item $\{\vp_j\}_{j\in J}$ tends to be equidistributed in the physical space.
\item $\{\vp_j\}_{j\in J}$ tends to be equidistributed in the phase space.
\end{enumerate}
\end{thm}

\begin{proof}
(ii)$\Rightarrow$(i). This is the same as in the proof of Theorem \ref{thm:QUE}.

(i)$\Rightarrow$(ii). Let $\nu$ be the semiclassical measure of a subsequence. By (i), $\Pi_*\nu=d\Vol_\Omega$. In particular, $\nu$ has total mass $1$. That is, there is no loss of mass at the boundary. Since $\nu$ is invariant, Theorem \ref{thm:billiard} gives $\nu=d\Vol_L$. Since every subsequence has a further subsequence with this limit, (ii) holds.
\end{proof}

In particular, $\Omega$ is QUE if and only if every sequence of normalized eigenfunctions with eigenvalues tending to infinity tends to be equidistributed in the physical space.

\begin{figure}[!ht]
\centering
\begin{tikzpicture}
\foreach \j in {0,...,9} \fill[gray!35] (-3,-1.9+0.4*\j) rectangle (3,-1.7+0.4*\j);
\draw[thick] (-3,2) -- (3,2);
\draw[thick] (-3,-2) -- (3,-2);
\draw[thick] (3,-2) arc (270:450:2);
\draw[thick] (-3,2) arc (90:270:2);
\draw[dashed] (-3,-2) -- (-3,2);
\draw[dashed] (3,-2) -- (3,2);
\draw[|-|] (-3,-2.4) -- (3,-2.4) node[midway,below] {$[-\ell,\ell]$};
\draw[|-|] (5.6,-2) -- (5.6,2) node[midway,right] {$[-1,1]$};
\node at (-4,0) {$\Omega$};
\end{tikzpicture}
\caption{A schematic bouncing ball mode. The eigenfunction concentrates in the rectangle $R$ and oscillates in the vertical direction. The shaded bands indicate where $|\vp|^2$ is large.}\label{fig:bb}
\end{figure}
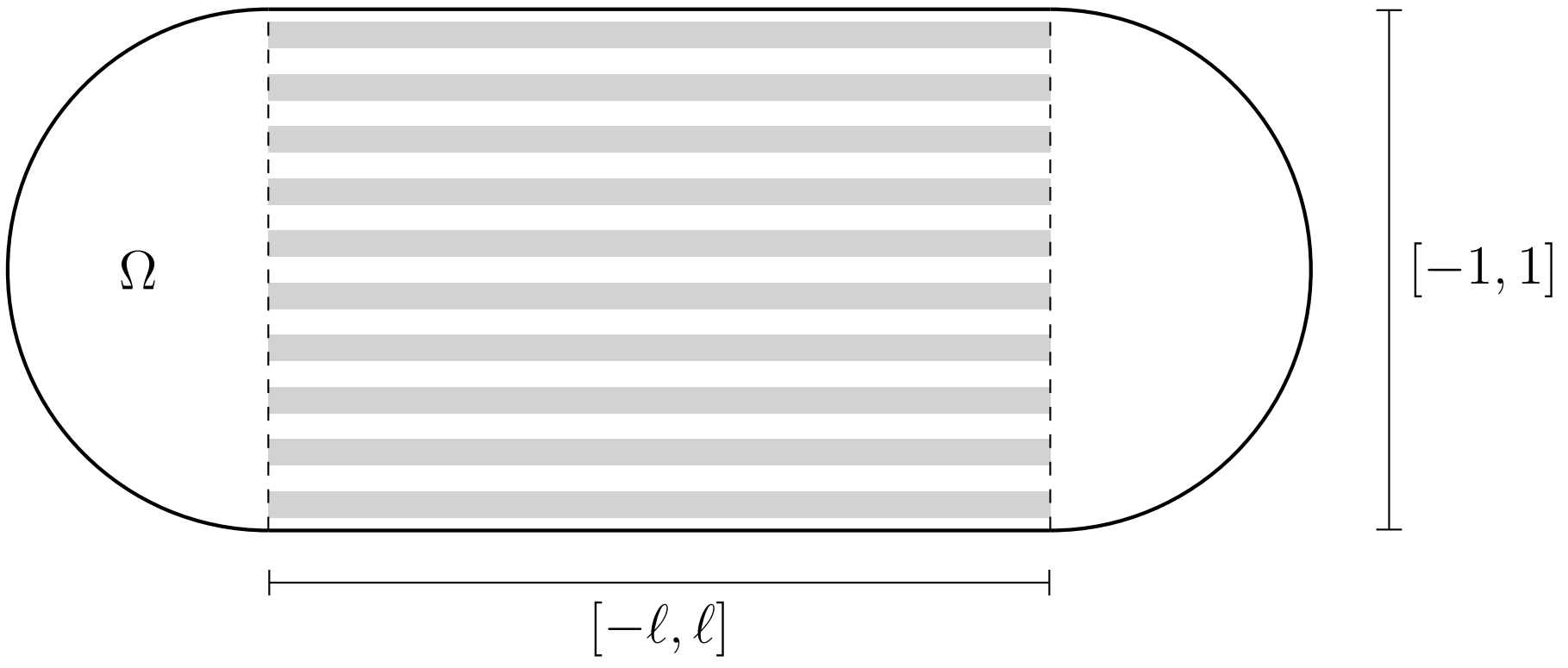

Heller \cite{Hel} observed numerically that some eigenfunctions of the stadium concentrate near unstable periodic orbits, and called this phenomenon ``\textit{scarring}''. Numerically, some eigenfunctions also concentrate in the rectangle $R$, near the family of bouncing ball orbits, which are neutral rather than unstable; see O'Connor--Heller \cite{OCH} and Figure \ref{fig:bb}. These observations prompted the quantum scarring conjecture: there are a subset $A\subset\Omega$ and a sequence of eigenfunctions such that
\[\int_A|\vp_j|^2\,d\Vol_\Omega\text{ does not converge to }\frac{|A|}{|\Omega|};\]
see Tao \cite[Conjecture 3.1]{T}. It is expected that $A$ can be taken as $R$. The eigenfunctions concentrating in $R$ are referred to as the bouncing ball modes.

Hassell \cite{Has} proved that $\Omega$ is not QUE for almost every $\ell\in[1,2]$. More precisely, for such $\ell$, some semiclassical measure gives positive mass to the family of bouncing ball trajectories in the phase space; see also Mangoubi--Weller Weiser \cite{MW}. Hassell's argument is variational in the parameter $\ell$, so it only covers almost every $\ell$. The result is expected to hold for every stadium, but this remains open.

By Theorem \ref{thm:QUEstadium}, a sequence of eigenfunctions which does not tend to be equidistributed in the phase space does not tend to be equidistributed in the physical space. The following theorem shows that this failure is detected on a disk. Combined with Hassell's theorem, it gives scarring in the physical space for almost every $\ell\in[1,2]$. This is a new proof of the positive answer to the quantum scarring conjecture for these $\ell$, which is due to Li--Liu \cite{LiLi}; see the discussion after the proof.
\begin{thm}[Physical-space scarring in the stadium]\label{thm:scar}
Let $\Omega$ be a Bunimovich stadium with the Dirichlet boundary condition. Suppose that a sequence of normalized eigenfunctions with eigenvalues tending to infinity does not tend to be equidistributed in the phase space. Then it has a subsequence $\{\vp_j\}_{j\in J}$, and there is an open disk $A$ with $\ol A\subset\Omega$, such that
\[\lim_{j\in J,j\to\infty}\int_A|\vp_j|^2\,d\Vol_\Omega\ne\frac{|A|}{|\Omega|}.\]
In particular, for almost every $\ell\in[1,2]$, the stadium with rectangular part $(-\ell,\ell)\times(-1,1)$ has such a sequence.
\end{thm}

\begin{proof}
Since the sequence does not tend to be equidistributed in the phase space, it has a subsequence $\{\vp_j\}_{j\in J}$ whose matrix elements converge to a semiclassical measure $\nu\ne d\Vol_L$. Then $\Pi_*\nu\ne d\Vol_\Omega$ by the proof of Theorem \ref{thm:QUEstadium}. That is,
\[\sigma=\Pi_*\nu-d\Vol_\Omega\]
is a non-zero finite signed measure on $\Omega$.

Write $D(x,r)$ for the open disk of radius $r$ centered at $x$. We claim that $\sigma(A)\ne0$ for some disk $A=D(x_0,r_0)$ with $\ol A\subset\Omega$. Suppose not. Let $f\in C_0(\Omega)$. For small $r>0$, denote by $f_r(x)$ the average of $f$ over $D(x,r)$. By Fubini's theorem,
\[\int_\Omega f_r\,d\sigma=\frac1{\pi r^2}\int_\Omega\int_\Omega\indicator_{D(y,r)}(x)f(y)\,dyd\sigma(x)=\frac1{\pi r^2}\int_\Omega f(y)\sigma\left(D(y,r)\right)dy=0.\]
Since $f_r\to f$ uniformly as $r\to0$, we obtain
\[\int_\Omega f\,d\sigma=0.\]
So $\sigma=0$, a contradiction.

The circles centered at $x_0$ are disjoint. So all but countably many of them have $\Pi_*\nu$-measure zero. Since $\sigma(D(x_0,r))\to\sigma(A)$ as $r\to r_0^-$, decreasing $r_0$ if necessary, we may assume that $\Pi_*\nu(\partial A)=0$. Since $|\vp_j|^2\,d\Vol_\Omega$ converges to $\Pi_*\nu$ tested against $C_0^\infty(\Omega)$, we obtain
\[\lim_{j\in J,j\to\infty}\int_A|\vp_j|^2\,d\Vol_\Omega=\Pi_*\nu(A)\ne\frac{|A|}{|\Omega|}.\]

For the last statement, let $\ell$ be as in Hassell's theorem. Then some semiclassical measure gives positive mass to the bouncing ball trajectories, which form a set of Liouville measure zero. So this semiclassical measure is not $d\Vol_L$, and the corresponding sequence of eigenfunctions does not tend to be equidistributed in the phase space. This completes the proof.
\end{proof}

That is, on every stadium, every failure of QUE is visible in the physical space, whatever the form of the semiclassical measure. However, Theorem \ref{thm:scar} does not locate the disk $A$. In particular, it does not show that the conclusion holds with $A$ replaced by the rectangle $R$, which contains the bouncing ball trajectories. Scarring is expected in $R$, but this is not known, even for almost every $\ell$, from either Theorem \ref{thm:scar} or the result of Li--Liu \cite{LiLi} below.

Li--Liu \cite[Theorem 1.1]{LiLi} proved that for almost every $\ell\in[1,2]$, there are normalized real Dirichlet eigenfunctions on $\Omega$ and a real function $f\in C^\infty(\ol\Omega)$ with zero mean such that
\[\int_\Omega f|\vp_j|^2\,d\Vol_\Omega\]
has a non-zero subsequential limit. Here, we have rescaled their stadium and eigenfunction normalization to match our conventions. Their proof takes Hassell's argument as an input and treats its cases separately, with an observable constructed in each case. Our argument recovers their theorem. Indeed, for almost every $\ell\in[1,2]$, we have $\Pi_*\nu\ne d\Vol_\Omega$ in the proof of Theorem \ref{thm:scar}. So there is a real function $f\in C_0^\infty(\Omega)$ such that
\begin{equation}\label{eq:scarob}
\int_\Omega f\,d\Pi_*\nu\ne\int_\Omega f\,d\Vol_\Omega.
\end{equation}
Subtracting from $f$ its mean with respect to $d\Vol_\Omega$, we obtain a function with zero mean. Since $\|\vp_j\|_{L^2(\Omega)}=1$, its integral against $|\vp_j|^2\,d\Vol_\Omega$ converges to a non-zero limit. Moreover, Theorem \ref{thm:scar} does not use the structure of Hassell's argument. It applies to every stadium and to every sequence of eigenfunctions which does not tend to be equidistributed in the phase space. On the other hand, our observable in \eqref{eq:scarob} is not explicit.

Scarring in the phase space is also known for the elliptical stadia in Figure \ref{fig:elliptic}. Fix $\ell_0>0$. For $b>0$, stretch the Bunimovich stadium with $\ell=\ell_0$ horizontally by the factor $b$. This gives the elliptical stadium with $\ell=b\ell_0$ and horizontal semi-axis $b$. Li--Zworski \cite[Theorem 3]{LiZw} proved that for almost every $b>0$, this stadium has a sequence of eigenfunctions whose semiclassical measures are supported on the closed bouncing ball family. Then the physical-space scarring follows directly, with $A$ close to the rectangle. This is a stronger form of scarring, and it is not covered by our results. Theorem \ref{thm:scar} detects the failure of equidistribution, but it gives no lower bound on its size. For the Bunimovich stadium, Burq--Zworski \cite{BZ} proved that the mass of every eigenfunction in any neighborhood of the closed half-disks is bounded below. So the eigenfunctions cannot concentrate in a compact subset of $R$.

\subsection*{Organization of the paper}
Section \ref{sec:surface} proves Theorem \ref{thm:surface} and extends it to finite-area hyperbolic orbifolds. Section \ref{sec:stadium} proves Theorem \ref{thm:billiard}. Section \ref{sec:higherproof} proves Theorem \ref{thm:entropy}. Appendix \ref{sec:toral} constructs invariant measures for hyperbolic toral automorphisms whose projection to the position is the Lebesgue measure.

%%%%%%%%%%%%%%%%%%%%%%%%%%%%%%%%%%%%%%%%%
\section{Surfaces}\label{sec:surface}
In this section, we prove Theorem \ref{thm:surface}. The proof uses the Fourier analysis in the fibers of $S\M$ and two facts about holomorphic functions on the transport twistor space. Throughout the section, let $\M$ be a connected complete oriented Riemannian surface of finite area.

\subsection{Fourier modes and the transport twistor space}\label{sec:twistor}
Denote by $Q$ the generator of the geodesic flow $a_t$. Write $\rho_\theta$ for the rotation by the angle $\theta\in\R$ in the fibers of $S\M$, which is defined by the orientation. Denote by $V$ the generator of $\rho_\theta$. In the coordinates $(x,\theta)$ given by a local oriented orthonormal frame, $V=\partial_\theta$, and the Liouville measure is
\[d\Vol_L=d\Vol_\M\otimes\frac{d\theta}{2\pi}.\]
Regard functions in $L^1(S\M)$ as distributions on $S\M$ by integration against $d\Vol_L$.

For $k\in\Z$, we say that a distribution $u$ on $S\M$ has weight $k$ if $Vu=iku$, that is,
\[u\circ\rho_\theta=e^{ik\theta}u\quad\text{for }\theta\in\R.\]
Every $u\in\Dc'(S\M)$ has the Fourier expansion
\[u=\sum_{k\in\Z}u_k,\quad\text{where }u_k\in\Dc'(S\M)\text{ has weight $k$ for }k\in\Z,\]
that is,
\[u_k=\frac1{2\pi}\int_0^{2\pi}e^{-ik\theta}u\circ\rho_\theta\,d\theta.\]
We say that $u$ is fiberwise holomorphic if $u_k=0$ for $k<0$.

Write $X_\perp=[Q,V]$. Following Guillemin--Kazhdan \cite{GK}, write
\[Q=\eta_++\eta_-,\quad\text{where }\eta_\pm=\frac12\left(Q\pm iX_\perp\right).\]
Then the first-order operators $\eta_\pm$ map distributions of weight $k$ to distributions of weight $k\pm1$. As a consequence, the weight-$k$ Fourier mode of $Qu$ is
\[(Qu)_k=\eta_+u_{k-1}+\eta_-u_{k+1}.\]
Thus, $Qu=0$ if and only if
\begin{equation}\label{eq:Fourier}
\eta_+u_{k-1}+\eta_-u_{k+1}=0\quad\text{for }k\in\Z.
\end{equation}
That is, the invariance under the geodesic flow is a recurrence relation between the Fourier modes.

Recall the closed unit ball bundle $Z$ in \eqref{eq:Z}. Write $Z^\circ$ for its interior. Denote by $\D$ the unit disk in $\C$. Write
\[\omega v=r\rho_\theta(v)\quad\text{for }\omega=re^{i\theta}\in\D\text{ and }v\in S_x\M.\]
The map $S\M\times\D\to Z^\circ$ defined by
\[(x,v,\omega)\mapsto(x,\omega v)\]
is onto, and its fibers are the orbits of the circle action
\[(x,v,\omega)\mapsto\left(\rho_\theta(x,v),e^{-i\theta}\omega\right)\quad\text{for }\theta\in\R.\]
So identify functions on $Z^\circ$ with functions on $S\M\times\D$ which are invariant under this action. For example, if $u_k$ has weight $k$ for each $k\in\N$, then
\[\left(u_k\circ\rho_\theta\right)\left(e^{-i\theta}\omega\right)^k=e^{ik\theta}u_ke^{-ik\theta}\omega^k=u_k\omega^k,\]
which implies that the series
\[\sum_{k\ge0}u_k(x,v)\omega^k\]
is invariant under the circle action.

Bohr--Paternain \cite{BP} introduced the transport twistor space. Here, we recall the description of its complex structure in Bohr--Lefeuvre--Paternain \cite[Section 3.1]{BLP}. A function on $Z^\circ$, regarded as a function $U$ on $S\M\times\D$, is holomorphic if and only if it satisfies
\begin{equation}\label{eq:CR}
\partial_{\ol\omega}U=0\quad\text{and}\quad\left(\omega^2\eta_++\eta_-\right)U=0.
\end{equation}
The first equation says that the restriction of a holomorphic function to each fiber $\{w\in T_x\M:|w|_x<1\}$ is holomorphic in $\omega$ in the usual sense. In \cite[Section 3.1]{BLP}, $\M$ may be any oriented Riemannian surface, compact or not. The complex structure degenerates at the boundary $\partial Z=S\M$.

The transport twistor space converts the flow-invariance equation into a holomorphicity condition. Indeed, suppose that
\[U=\sum_{k\ge0}u_k\omega^k,\quad\text{where $u_k$ has weight $k$ for }k\in\N.\]
Shifting the index in the first sum, we obtain
\[\left(\omega^2\eta_++\eta_-\right)U=\sum_{k\ge0}\left(\eta_+u_{k-2}+\eta_-u_k\right)\omega^k,\]
where $u_{-1}=u_{-2}=0$. By \eqref{eq:Fourier}, these coefficients vanish if and only if the fiberwise holomorphic distribution with Fourier modes $u_k$ is invariant under the geodesic flow.

\subsection{Holomorphic functions on the transport twistor space}
The following two lemmas only use the local equations \eqref{eq:CR}, so they apply to compact and non-compact surfaces alike. The first lemma constructs a holomorphic function from an invariant distribution. The second lemma uses the ergodicity to show that every bounded holomorphic function is constant.

\begin{lemma}[Holomorphic extension]\label{lemma:extension}
Let $u\in\Dc'(S\M)$ be fiberwise holomorphic, that is,
\[u=\sum_{k\ge0}u_k,\quad\text{where }u_k\in\Dc'(S\M)\text{ has weight $k$ for }k\in\N.\]
Suppose that $Qu=0$ and
\[\sup_{k\ge0}\left\|u_k\right\|_{L^\infty(S\M)}<\infty.\]
Then
\[U(x,v,\omega)=\sum_{k\ge0}u_k(x,v)\omega^k\quad\text{for }(x,v)\in S\M\text{ and }\omega\in\D\]
defines a holomorphic function on $Z^\circ$.
\end{lemma}

\begin{proof}
Fix $0<r<1$. Since the coefficients $u_k$ are uniformly bounded, the series converges absolutely in $L^\infty(S\M\times\{|\omega|\le r\})$. Each term is invariant under the circle action, as shown in Section \ref{sec:twistor}. So $U$ is a locally bounded function on $Z^\circ$, and in particular, a distribution on $Z^\circ$.

Since $U$ is a power series in $\omega$ whose coefficients do not depend on $\omega$, we have $\partial_{\ol\omega}U=0$. The operators $\eta_\pm$ differentiate only in $(x,v)$, and the series converges locally uniformly in $\omega$ with values in $L^\infty(S\M)$. So $\eta_\pm$ can be applied term by term in distributions. As in Section \ref{sec:twistor},
\[\left(\omega^2\eta_++\eta_-\right)U=\sum_{k\ge0}\left(\eta_+u_{k-2}+\eta_-u_k\right)\omega^k,\quad\text{where }u_{-1}=u_{-2}=0.\]
By \eqref{eq:Fourier} with $k-1$ in place of $k$, each coefficient $\eta_+u_{k-2}+\eta_-u_k$ is the weight-$(k-1)$ Fourier mode of $Qu=0$, so it vanishes. Hence,
\[\left(\omega^2\eta_++\eta_-\right)U=0.\]
That is, $U$ satisfies the Cauchy--Riemann equations \eqref{eq:CR} in the sense of distributions on the complex surface $Z^\circ$. By the elliptic regularity of the Cauchy--Riemann equations, $U$ is a holomorphic function. This proves the lemma.
\end{proof}

\begin{lemma}[Bounded holomorphic functions]\label{lemma:bounded}
Suppose that the geodesic flow on $S\M$ is ergodic with respect to $d\Vol_L$. Then every bounded holomorphic function on $Z^\circ$ is constant.
\end{lemma}

The idea is to take the boundary value of a bounded holomorphic function on $Z^\circ$. The Cauchy--Riemann equations imply that the boundary value is invariant under the geodesic flow. The ergodicity forces it to be constant. Since the boundary value determines the original holomorphic function, that function is constant.

\begin{proof}
Let $F$ be a bounded holomorphic function on $Z^\circ$, regarded as a function on $S\M\times\D$. By \eqref{eq:CR}, $F$ is holomorphic in $\omega$. Write
\[F(x,v,\omega)=\sum_{k\ge0}b_k(x,v)\omega^k\quad\text{for }(x,v)\in S\M\text{ and }\omega\in\D.\]
Here,
\[b_k=\frac1{k!}\partial_\omega^kF\big|_{\omega=0}.\]
Then $b_k\in C^\infty(S\M)$. For fixed $(x,v)$, the function $\omega\mapsto F(x,v,\omega)$ is bounded and holomorphic on $\D$. So the Cauchy estimates give
\[\left\|b_k\right\|_{L^\infty(S\M)}\le\|F\|_{L^\infty(Z^\circ)}.\]
Since $F$ is invariant under the circle action, comparing the coefficients of $\omega^k$ in
\[F\left(\rho_\theta(x,v),e^{-i\theta}\omega\right)=F(x,v,\omega)\]
shows that
\[b_k\circ\rho_\theta=e^{ik\theta}b_k,\]
that is, $b_k$ has weight $k$.

Now take the boundary value of $F$. For $0<r<1$, write
\[F_r(x,v)=F(x,v,r)=\sum_{k\ge0}r^kb_k(x,v).\]
Then
\[\|F_r\|_{L^\infty(S\M)}\le\|F\|_{L^\infty(Z^\circ)}.\]
By the second equation in \eqref{eq:CR} at $\omega=r$, we have
\[\eta_-F_r=-r^2\eta_+F_r.\]
So
\[QF_r=\left(1-r^2\right)\eta_+F_r.\]
That is, the derivative of $F_r$ along the geodesic flow carries the factor $1-r^2$, which vanishes at the boundary of the fibers. Since $Q$ and $X_\perp$ preserve $d\Vol_L$, integrating by parts twice gives
\[\int_{S\M}F_rQ\phi\,d\Vol_L=\left(1-r^2\right)\int_{S\M}F_r\eta_+\phi\,d\Vol_L\quad\text{for }\phi\in C_0^\infty(S\M).\]
Hence,
\[\left|\int_{S\M}F_rQ\phi\,d\Vol_L\right|\le\left(1-r^2\right)\|F\|_{L^\infty(Z^\circ)}\left\|\eta_+\phi\right\|_{L^1(S\M)}.\]
Since the functions $F_r$ are uniformly bounded in $L^\infty(S\M)$, choose $r_j\to1$ such that $F_{r_j}$ converges weak-* to some $f\in L^\infty(S\M)$. Letting $r=r_j\to1$, we obtain
\[\int_{S\M}fQ\phi\,d\Vol_L=0\quad\text{for }\phi\in C_0^\infty(S\M).\]
That is, $Qf=0$ in $\Dc'(S\M)$.

Since $Qf=0$ and the geodesic flow preserves $d\Vol_L$, we have $f\circ a_t=f$ almost everywhere for $t\in\R$. By the ergodicity of the geodesic flow, $f$ is constant almost everywhere. The projection $u\mapsto u_k$ to the weight-$k$ Fourier mode is continuous in the weak-* topology, and the weight-$k$ Fourier mode of $F_{r_j}$ is $r_j^kb_k$. So the weight-$k$ Fourier mode of $f$ is $b_k$ for $k\ge0$. Since $f$ is constant, $b_k=0$ for $k\ge1$, and $F=b_0$ is constant. This proves the lemma.
\end{proof}

\subsection{Proof of Theorem \ref{thm:surface}}\label{sec:surfaceproof}
\begin{figure}[!ht]
\centering
\begin{tikzpicture}
\draw[thick] (0,0) circle (3);
\draw[dashed] (0,0) -- (150:3) node[midway,below left] {$1$};
\draw[dashed] (0,0) -- (1.2,.9) node[midway,below right] {$|\omega|$};
\draw[thick] (150:3) -- (1.2,.9) node[midway,above] {$|\zeta-\omega|$};
\fill (0,0) circle (.06) node[below] {$0$};
\fill (1.2,.9) circle (.06) node[right] {$\omega$};
\fill (150:3) circle (.06) node[left] {$\zeta$};
\node at (3,-2.6) {$S_x\M$};
\end{tikzpicture}
\caption{The fiber of $Z^\circ$ over $x\in\M$ is the unit disk in $T_x\M$, identified with $\D$ by $\omega v\mapsto\omega$. Its boundary $S_x\M$ carries the conditional measure $\nu_x$. The restriction of $H$ to this fiber is the Carath\'eodory transform of $\nu_x$. For $\zeta\in S_x\M$ and $\omega\in\D$, the triangle inequality gives $|\zeta-\omega|\le1+|\omega|$.}\label{fig:fiber}
\end{figure}
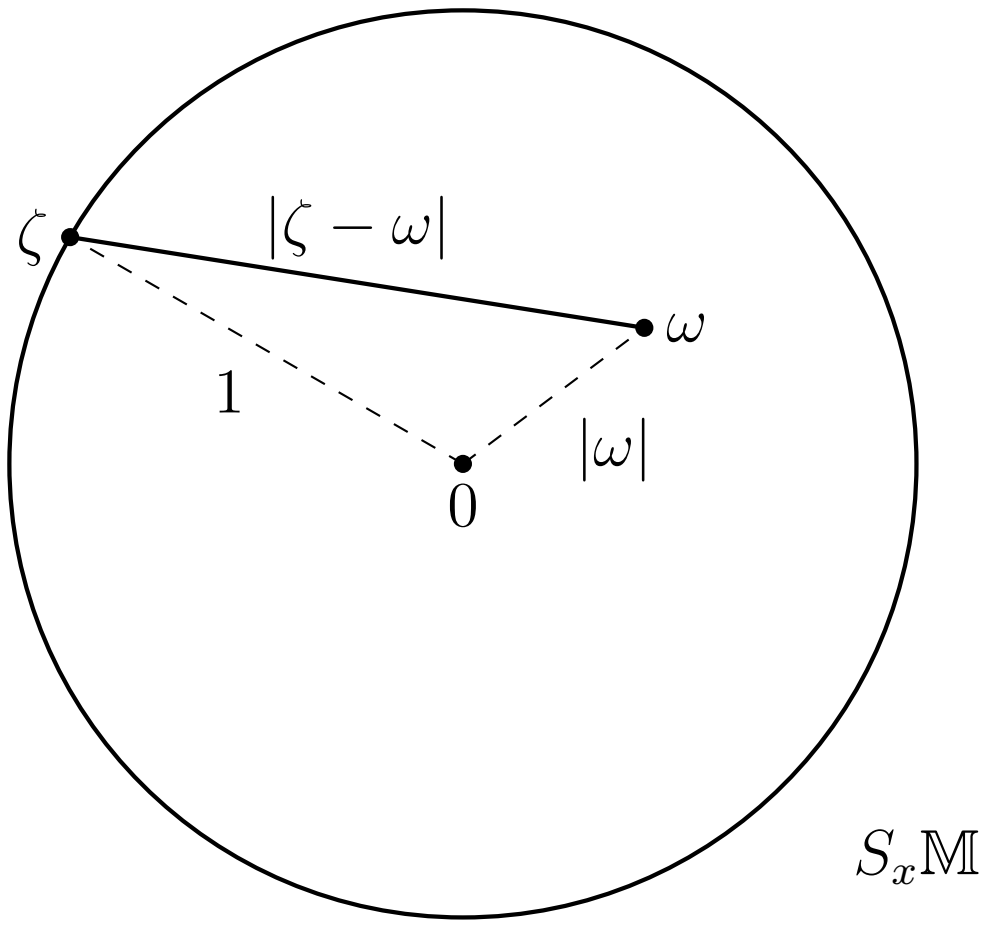

If $\mu$ has a density in $L^2(S\M)$, then the proof is short. Indeed, the fiberwise holomorphic distribution $u$ in Step 2 below is then an invariant function in $L^2(S\M)$, so it is constant by the ergodicity. In general, $\mu$ may be singular. In this case, $u$ is only an invariant distribution, and the ergodicity does not apply. The transport twistor space overcomes this difficulty. The idea is to encode the conditional measures of $\mu$ on the fibers of $S\M$ in one holomorphic function on $Z^\circ$. The fiber of $Z^\circ$ over $x\in\M$ is a disk with boundary circle $S_x\M$; see Figure \ref{fig:fiber}. A probability measure on a circle is encoded by its Carath\'eodory transform, which is a holomorphic function on the disk with positive real part. By Section \ref{sec:twistor}, these functions on the fibers form a holomorphic function on $Z^\circ$ if the corresponding fiberwise holomorphic distribution is invariant. Then Lemma \ref{lemma:bounded} forces this holomorphic function to be trivial. The proof has five steps.

\textit{Step 1: Removal of the odd Fourier modes.} By Section \ref{sec:twistor}, holomorphic functions on $Z^\circ$ correspond to fiberwise holomorphic invariant distributions. So we aim to keep only the Fourier modes of $\mu$ of non-negative weight. This truncation must preserve the invariance. By \eqref{eq:Fourier}, the equation of weight $k$ only involves the Fourier modes of weights $k-1$ and $k+1$. For the truncation, the equations of weights $k\ge2$ are the same as those for $\mu$, and the equations of weights $k\le-2$ hold trivially. However, the truncation breaks the equation of weight $0$, which couples the Fourier modes of weights $-1$ and $1$. Indeed, denote by $Y$ the mean velocity field of $\mu$. That is, $Y(x)$ is the mean of $v$ with respect to the conditional measure of $\mu$ on $S_x\M$. For $\mu$, the equation of weight $0$ says that $\Div Y=0$. For the truncation, it says that $\Div Y=0$ and $\curl Y=0$. Here, we used that the Fourier mode $\mu_1$ of $\mu$ of weight $1$ satisfies
\[\eta_-\mu_1=\frac12\left(\Div Y-i\curl Y\right).\]
The invariance does not imply $\curl Y=0$. The equations of even weight only involve the odd Fourier modes. So if the odd Fourier modes vanish, then all these equations hold trivially, and the problem disappears. The equations of weights $\pm1$ involve the zeroth Fourier mode, which is constant by the uniform projection; see Step 2. Therefore, first remove the odd Fourier modes.

Denote by $R(x,v)=(x,-v)$ the velocity reversal. Since $R\circ a_t=a_{-t}\circ R$, the measure $R_*\mu$ satisfies \eqref{eq:IP}. Let
\[\nu=\frac12\left(\mu+R_*\mu\right).\]
Then $\nu$ is invariant under $a_t$, satisfies $\Pi_*\nu=d\Vol_\M$, and $R_*\nu=\nu$. Since $\Pi_*\nu=d\Vol_\M$, disintegrate
\[d\nu(x,v)=d\nu_x(v)\,d\Vol_\M(x),\]
where $\nu_x$ is a probability measure on $S_x\M$ for almost every $x\in\M$. Regard $\nu$ as a distribution on $S\M$. Write
\[\nu=\sum_{k\in\Z}\nu_k.\]
For almost every $x\in\M$ and $v\in S_x\M$, identify $S_x\M$ with the unit circle $\Sb$ in $\C$ by $\rho_\theta(v)\mapsto e^{i\theta}$. The Fourier modes of $\nu$ are the functions
\[\nu_k(x,v)=\int_{\Sb}\zeta^{-k}\,d\nu_x(\zeta).\]
Since $\nu_x$ is a probability measure for almost every $x$, we have $\nu_0=1$, and $|\nu_k|\le1$ almost everywhere for $k\in\Z$. Since $R=\rho_\pi$ and $R_*\nu=\nu$, we have $\nu_k=(-1)^k\nu_k$. That is, $\nu_k=0$ for odd $k$. Since $\nu$ is real, $\nu_{-k}=\ol{\nu_k}$ for $k\in\Z$. So $\nu$ is determined by its Fourier modes of non-negative weight, and the truncation loses no information.

\textit{Step 2: A fiberwise holomorphic invariant distribution.} Let
\[u=1+2\sum_{k\ge1}\nu_{2k}.\]
The coefficient $2$ is explained in Step 3. Then $u$ is fiberwise holomorphic, and $|u_k|\le2$ almost everywhere for $k\in\N$. That is, $u$ satisfies the bound in Lemma \ref{lemma:extension}. We claim that $Qu=0$. Since $u$ has only non-negative even Fourier modes, all even-weight equations and all equations of weight at most $-3$ hold automatically. Since $u_0=1$, the equation of weight $-1$ reads $\eta_-1=0$, which holds. Since $(a_t)_*\nu=\nu$, we have $Q\nu=0$. Since $\nu_0=1$, the equation of weight $1$ in \eqref{eq:Fourier} for $\nu$ gives
\[\eta_-\nu_2=-\eta_+\nu_0=0,\]
and hence
\[\eta_+u_0+\eta_-u_2=0.\]
For $k\ge1$, the remaining equations are
\[\eta_+u_{2k}+\eta_-u_{2k+2}=2\left(\eta_+\nu_{2k}+\eta_-\nu_{2k+2}\right)=0.\]
Thus, $Qu=0$.

This is where the uniform projection enters. Suppose instead that $\Pi_*\nu$ has a non-constant smooth density $\rho$ with respect to $d\Vol_\M$. Then $\nu_0=\rho$. The truncation
\[\rho+2\sum_{k\ge1}\nu_{2k}\]
satisfies the equations of weights $-1$ and $1$ only if $\eta_-\rho=\eta_+\rho=0$, that is, only if $\rho$ is constant.

\textit{Step 3: The Carath\'eodory transform.} By Lemma \ref{lemma:extension}, we obtain the holomorphic function
\[H(x,v,\omega)=1+2\sum_{k\ge1}\nu_{2k}(x,v)\omega^{2k}\quad\text{on }Z^\circ.\]
Now identify $H$ on each fiber. Recall that the Carath\'eodory transform of a probability measure $\sigma$ on the unit circle $\Sb\subset\C$ is the holomorphic function
\[\omega\mapsto\int_{\Sb}\frac{\zeta+\omega}{\zeta-\omega}\,d\sigma(\zeta)\quad\text{for }\omega\in\D;\]
see Simon \cite[Section 1.3]{Sim}. It is also called the Herglotz transform. Its real part is the Poisson integral of $\sigma$, and it takes the value $1$ at $\omega=0$. Since
\[\frac{\zeta+\omega}{\zeta-\omega}=1+2\sum_{k\ge1}\zeta^{-k}\omega^k\quad\text{for }|\zeta|=1\text{ and }|\omega|<1,\]
and the odd Fourier modes of $\nu$ vanish, we obtain
\[H(x,v,\omega)=\int_{\Sb}\frac{\zeta+\omega}{\zeta-\omega}\,d\nu_x(\zeta)\quad\text{for }\omega\in\D.\]
That is, the restriction of $H$ to the fiber of $Z^\circ$ over $x$ is the Carath\'eodory transform of $\nu_x$; see Figure \ref{fig:fiber}. This explains the coefficient $2$ in the definition of $u$. It converts the positivity of $\nu_x$ into the positivity of $\Rer H$.

Since $\nu_x$ is a probability measure,
\[\Rer H(x,v,\omega)=\int_{\Sb}\frac{1-|\omega|^2}{|\zeta-\omega|^2}\,d\nu_x(\zeta)\ge\frac{1-|\omega|}{1+|\omega|}\quad\text{for }\omega\in\D.\]
Here, we used the triangle inequality $|\zeta-\omega|\le1+|\omega|$; see Figure \ref{fig:fiber}. Since $H$ is continuous on $Z^\circ$, this inequality holds for every $x\in\M$. In particular, $\Rer H>0$ on $Z^\circ$.

\textit{Step 4: The Cayley transform and the ergodicity.} Lemma \ref{lemma:bounded} applies only to bounded holomorphic functions. In general, $H$ is unbounded near $\partial Z$. For example, if $\nu_x$ has an atom at $\zeta_0\in\Sb$, then $\Rer H(x,v,r\zeta_0)\to\infty$ as $r\to1^-$. Use the positivity of $\Rer H$ to obtain a bounded holomorphic function. The Cayley transform
\[s\mapsto\frac{s-1}{s+1}\]
maps the right half-plane $\{\Rer s>0\}$ onto $\D$. Since $\Rer H>0$,
\[F=\frac{H-1}{H+1}\]
is holomorphic on $Z^\circ$ and satisfies $|F|<1$. By Lemma \ref{lemma:bounded}, $F$ is constant. The normalization on the zero section determines the constant. Indeed, since the zeroth Fourier mode of $u$ is $1$, we have $H(x,v,0)=1$, so $F=0$. That is, $H=1$. Hence, for $k\ge1$,
\[\nu_{2k}=0\quad\text{and}\quad\nu_{-2k}=\ol{\nu_{2k}}=0.\]
So $\nu=d\Vol_L$.

\textit{Step 5: Extremality.} Step 1 replaced $\mu$ by $\nu$ and discarded the odd Fourier modes of $\mu$. Finally, recover $\mu$. Notice that
\[d\Vol_L=\frac12\left(\mu+R_*\mu\right),\]
where $\mu$ and $R_*\mu$ are invariant probability measures. Since $d\Vol_L$ is ergodic, it is an extreme point of the convex set of invariant probability measures. So $\mu=R_*\mu=d\Vol_L$. This completes the proof of Theorem \ref{thm:surface}.

\begin{rmk}[The role of the ergodicity]
The ergodicity is used twice, in Steps 4 and 5. Steps 1--3 and the construction of $F$ apply to every probability measure satisfying \eqref{eq:IP}. The flat torus $\T^2=\R^2/\Z^2$ in the Example in Section \ref{sec:intro} shows what fails without the ergodicity. Identify each tangent space of $\T^2$ with $\C$. Then $S\T^2=\T^2\times\Sb$, and the geodesic flow is
\[a_t(x,v)=(x+tv,v).\]
It is not ergodic, since every function of $v$ alone is invariant.

Fix $v_0\in\Sb$. Let $\mu=dx\otimes\delta_{v_0}$. Then
\[\nu=dx\otimes\frac12\left(\delta_{v_0}+\delta_{-v_0}\right).\]
In the identification of $S_x\T^2$ with $\Sb$ in Step 1, the directions $\pm v_0$ correspond to $\pm v_0/v$. So
\[\nu_{2k}(x,v)=\left(\frac v{v_0}\right)^{2k}\quad\text{for }k\in\N.\]
Write $w=\omega v$. Steps 2 and 3 give
\[H(x,v,\omega)=1+2\sum_{k\ge1}\left(\frac w{v_0}\right)^{2k}=\frac{v_0^2+w^2}{v_0^2-w^2}.\]
Hence,
\[F(x,v,\omega)=\left(\frac w{v_0}\right)^2.\]
This is a non-constant bounded holomorphic function on $Z^\circ=\T^2\times\D$. Its boundary value in the proof of Lemma \ref{lemma:bounded} is
\[f(x,v)=\left(\frac v{v_0}\right)^2.\]
The function $f$ is invariant under $a_t$, but it is not constant. So the conclusion of Lemma \ref{lemma:bounded} fails, and $\nu\ne d\Vol_L$. In general, $F=0$ if and only if $\nu=d\Vol_L$. So $F$ measures the deviation of $\nu$ from $d\Vol_L$. In Step 4, the ergodicity rules out non-constant invariant boundary values such as $f$, and hence such deviations.

Step 5 also requires the ergodicity. Write $v=e^{i\theta}$. Let
\[d\mu=\left(1+\cos\theta\right)\frac{dxd\theta}{2\pi}.\]
Then $\mu$ is invariant under $a_t$, and $\Pi_*\mu=dx$. Since
\[R_*\mu=\frac{1-\cos\theta}{2\pi}\,dxd\theta,\]
we have $\nu=d\Vol_L$. So Steps 1--4 hold with $H=1$, but $\mu\ne d\Vol_L$. That is, the ergodicity in Step 5 recovers the odd Fourier modes of $\mu$ discarded in Step 1.
\end{rmk}

\subsection{Hyperbolic orbifolds}\label{sec:orbifold}
Theorem \ref{thm:surface} extends to finite-area hyperbolic orbifolds by passage to finite covers. Let $\Gamma\subset\PSL(2,\R)$ be a lattice. That is, $\Gamma$ is a discrete subgroup with $\Vol(\Gamma\bs\Hb)<\infty$. Write $\M=\Gamma\bs\Hb$. Identify $S\M$ with $\Gamma\bs\PSL(2,\R)$. If $\Gamma$ is torsion-free, then $\M$ is a complete hyperbolic surface of finite area, so Theorem \ref{thm:surface} applies to $\M$ by Hopf \cite{Ho1, Ho2}. If $\Gamma$ has torsion, then $\M$ is an orbifold. For example, this is the case for $\Gamma=\PSL(2,\Z)$. Still write $\Pi:S\M\to\M$ for the natural projection. Write $d\Vol_\M$ and $d\Vol_L$ for the normalized measures induced by the hyperbolic area and the Haar measure, respectively.

\begin{prop}[Lattices with torsion]\label{prop:lattice}
Let $\Gamma\subset\PSL(2,\R)$ be a lattice. If a probability measure $\mu$ on $S\M=\Gamma\bs\PSL(2,\R)$ satisfies \eqref{eq:IP}, then $\mu=d\Vol_L$.
\end{prop}

\begin{proof}
Lift $\mu$ to a finite cover by averaging over the fibers of the covering map. By Selberg's lemma, $\Gamma$ has a torsion-free subgroup $\tilde\Gamma$ of finite index $d$; see \cite[Lemma 8]{Sel}. For $\Gamma=\PSL(2,\Z)$, we may take the image of the principal congruence subgroup of level $2$. Write $\tilde\M=\tilde\Gamma\bs\Hb$. Denote by $p:S\tilde\M=\tilde\Gamma\bs\PSL(2,\R)\to S\M$ the $d$-sheeted covering. Define a probability measure $\tilde\mu$ on $S\tilde\M$ by
\[\int_{S\tilde\M}f\,d\tilde\mu=\int_{S\M}\frac1d\sum_{p(\tilde z)=z}f(\tilde z)\,d\mu(z)\quad\text{for }f\in C_0(S\tilde\M).\]
Since $p$ commutes with the geodesic flow, $\tilde\mu$ is invariant under the geodesic flow. The covering $\tilde\M\to\M$ has degree $d$ and is unramified outside the finitely many elliptic points, which have measure zero. So the projection of $\tilde\mu$ to $\tilde\M$ is the normalized hyperbolic area. By Theorem \ref{thm:surface}, $\tilde\mu$ is the Liouville measure on $S\tilde\M$. Hence,
\[\mu=p_*\tilde\mu=d\Vol_L.\]
This completes the proof.
\end{proof}
The same argument applies to compact non-orientable negatively curved surfaces, using the orientation double cover.

%%%%%%%%%%%%%%%%%%%%%%%%%%%%%%%%%%%%%%%%%
\section{Billiards}\label{sec:stadium}
In this section, we prove Theorem \ref{thm:billiard}. We follow the five steps in Section \ref{sec:surfaceproof} and point out the differences. Throughout the section, let $\Omega$ be as in Theorem \ref{thm:billiard}. Identify $x=(x_1,x_2)\in\R^2$ with $z=x_1+ix_2\in\C$. By an oriented chord of $\Omega$, we mean a maximal open line segment in $\Omega$ with a direction. Its endpoints lie on $\partial\Omega$. Between two consecutive reflections, a trajectory of the billiard flow runs along an oriented chord.

Compared with surfaces, there are three differences.

First, the billiard flow jumps at the reflections. So its invariance cannot be tested against arbitrary smooth functions. In Section \ref{sec:weak}, we test it against functions which take the same value on the incoming and the outgoing directions.

Second, the table is flat. So the Fourier modes, the operators $\eta_\pm$, and the transport twistor space are explicit. In Sections \ref{sec:flat} and \ref{sec:interior}, we carry out Steps 1--3 explicitly.

Third, the ergodicity is used differently in Step 4. In Lemma \ref{lemma:bounded}, the boundary value $f$ of $F$ is invariant under the geodesic flow. On a billiard table, $f$ is constant along the chords, but the reflections reverse the orientation of the fibers. In Sections \ref{sec:difficulty} and \ref{sec:Ref}, we show that $\Rer f$ is still invariant under the billiard flow, and that this suffices.

\begin{figure}[!ht]
\centering
\begin{tikzpicture}
\draw[thick] (-20:3) arc (-20:70:3);
\fill (2.719,1.268) circle (.05) node[above left] {$x$};
\draw[dotted,thick] (2.085,2.627) -- (3.353,-.091);
\draw[->,thick] (.749,1.616) -- (2.719,1.268);
\draw[->,thick] (2.719,1.268) -- (1.72,-.464);
\draw[->,thick] (2.719,1.268) -- (3.535,1.648) node[right] {$e^{i\psi(x)}$};
\node at (1.4,1.75) {$\theta$};
\node at (1.65,.2) {$\tau_x(\theta)$};
\node at (.4,-.3) {$\Omega$};
\end{tikzpicture}
\caption{The reflection at $x\in\partial\Omega$. The incoming direction $\theta$ and the outgoing direction $\tau_x(\theta)$ are symmetric with respect to the tangent line (dotted), and $e^{i\psi(x)}$ is the outward unit normal.}\label{fig:reflection}
\end{figure}
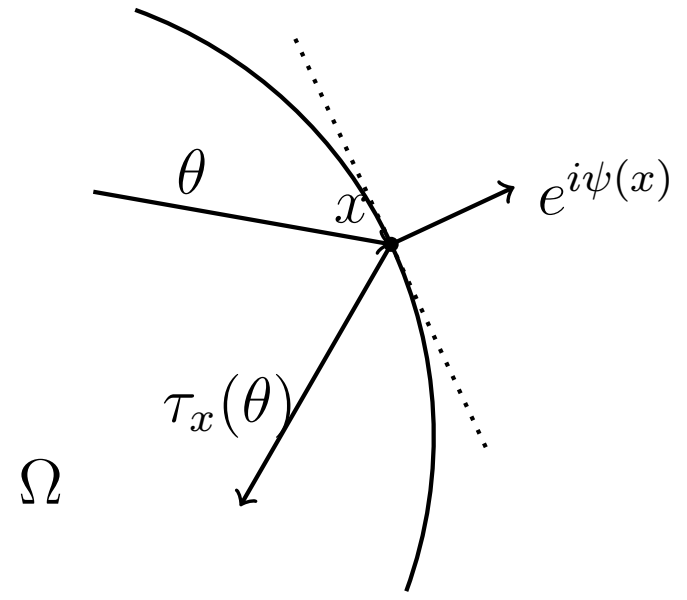

For $x\in\partial\Omega\setminus\Sigma$, denote by $\psi(x)$ the angle of the outward unit normal. So the outward unit normal is $e^{i\psi(x)}$. Let
\[\tau_x(\theta)=2\psi(x)+\pi-\theta\]
be the reflection of directions across the tangent line at $x$. Notice that
\begin{equation}\label{eq:cos}
\cos\left(\tau_x(\theta)-\psi(x)\right)=-\cos\left(\theta-\psi(x)\right).
\end{equation}
That is, the reflection reverses the normal component of the direction; see Figure \ref{fig:reflection}.

\subsection{Invariance in the weak form}\label{sec:weak}
On a surface, the invariance of a measure $\mu$ under the geodesic flow is equivalent to $Q\mu=0$ in $\Dc'(S\M)$. That is,
\[\int Q\phi\,d\mu=0\quad\text{for }\phi\in C_0^\infty(S\M).\]
On a billiard table, the flow jumps at the reflections. For a smooth function $\phi$ on $\ol\Omega\times\Sb$, the function $t\mapsto\phi(a_t(x,\theta))$ is discontinuous at a reflection, unless $\phi$ takes the same value on the incoming and the outgoing directions. So we only use test functions with this property.

Let $C_\tau$ be the space of real functions $\phi\in C^1(\ol\Omega\times\Sb)$ which vanish near $\Sigma\times\Sb$ and satisfy $\phi(x,\tau_x(\theta))=\phi(x,\theta)$ for $x\in\partial\Omega\setminus\Sigma$ and $\theta\in\Sb$. Write
\[Q=\cos\theta\partial_{x_1}+\sin\theta\partial_{x_2}\]
for the generator of the free motion. The following lemma replaces the equation $Q\mu=0$ on surfaces. Its part (ii) replaces the last step of the proof of Lemma \ref{lemma:bounded}, where $Qf=0$ implies that $f$ is invariant.

\begin{lemma}[Invariance in the weak form]\label{lemma:weak}\leavevmode
\begin{enumerate}[(i)]
\item Let $\mu$ be an invariant probability measure on $\ol\Omega\times\Sb$. Then
\[\int Q\phi\,d\mu=0\quad\text{for }\phi\in C_\tau.\]
\item Let $g\in L^\infty(\Omega\times\Sb)$. Suppose that
\[\int gQ\phi\,d\Vol_L=0\quad\text{for }\phi\in C_\tau.\]
Then $g\circ a_t=g$ almost everywhere for $t\in\R$.
\end{enumerate}
\end{lemma}

\begin{proof}\leavevmode
\begin{enumerate}[(i)]
\item Since $\mu$ is invariant, $a_t$ is defined $\mu$-almost everywhere. So $\mu$-almost every trajectory does not hit the corners and has finitely many reflections in a bounded time interval. For such a trajectory, the function $t\mapsto\phi(a_t(x,\theta))$ is continuous at the reflections. Indeed, at a reflection at $x\in\partial\Omega\setminus\Sigma$, this follows from $\phi(x,\tau_x(\theta))=\phi(x,\theta)$. Near a reflection at a point of $\Sigma$, this function vanishes, since $\phi$ vanishes near $\Sigma\times\Sb$. So it is Lipschitz, with derivative $(Q\phi)(a_t(x,\theta))$ for almost every $t$. By the dominated convergence theorem and the invariance of $\mu$,
\[0=\frac{d}{dt}\int\phi\circ a_t\,d\mu=\int(Q\phi)\circ a_t\,d\mu=\int Q\phi\,d\mu.\]
This proves (i).
\item Taking $\phi$ supported in $\Omega\times\Sb$, we obtain $Qg=0$ in $\Dc'(\Omega\times\Sb)$. So $g$ is constant along almost every oriented chord of $\Omega$. Write $g^\partial(x,\theta)$ for the value of $g$ on the chord through $x\in\partial\Omega\setminus\Sigma$ in the direction $\theta$. In the coordinates given by the arc length on $\partial\Omega$, the direction, and the time along the chords, we have
\[dxd\theta=\left|\cos\left(\theta-\psi(x)\right)\right|\,dsd\theta dt.\]
Integrating $Q\phi$ along each chord, we obtain
\[0=\int gQ\phi\,dxd\theta=\int_{\partial\Omega}\int_{\Sb}g^\partial(x,\theta)\phi(x,\theta)\cos\left(\theta-\psi(x)\right)\,d\theta ds(x).\]
Let $x_0\in\partial\Omega\setminus\Sigma$. Extend $\psi$ to a $C^1$ function near $x_0$ which is constant along the normals. For $\phi_1\in C_0^\infty(\R^2)$ supported near $x_0$ and $\phi_2\in C^\infty(\Sb)$, take
\[\phi(x,\theta)=\phi_1(x)\left(\phi_2(\theta)+\phi_2\left(\tau_x(\theta)\right)\right).\]
Then $\phi\in C_\tau$. By \eqref{eq:cos}, we obtain
\[g^\partial\left(x,\tau_x(\theta)\right)=g^\partial(x,\theta)\quad\text{for almost every }(x,\theta).\]
That is, $g$ takes the same value on the incoming and the outgoing chords at almost every reflection. So $g$ is invariant under the billiard flow. This proves (ii).
\end{enumerate}
\end{proof}

\subsection{Fourier modes and the transport twistor space of a flat table}\label{sec:flat}
On a flat table, the objects in Section \ref{sec:twistor} are explicit. Write $v=e^{i\theta}$. The rotation $\rho_\theta$ shifts the angle, so $V=\partial_\theta$. A distribution of weight $k$ on $\Omega\times\Sb$ has the form $g(x)e^{ik\theta}$ with $g\in\Dc'(\Omega)$. The generator of the free motion is
\[Q=e^{i\theta}\partial_z+e^{-i\theta}\partial_{\bar z}.\]
So the operators in Section \ref{sec:twistor} are
\[\eta_+=e^{i\theta}\partial_z\quad\text{and}\quad\eta_-=e^{-i\theta}\partial_{\bar z}.\]
They act on the Fourier modes by
\[\eta_+\left(ge^{ik\theta}\right)=\left(\partial_zg\right)e^{i(k+1)\theta}\quad\text{and}\quad\eta_-\left(ge^{ik\theta}\right)=\left(\partial_{\bar z}g\right)e^{i(k-1)\theta}.\]
Hence, for
\[u=\sum_{k\in\Z}g_k(x)e^{ik\theta},\]
the recurrence \eqref{eq:Fourier} reads
\begin{equation}\label{eq:Fourierflat}
\partial_zg_{k-1}+\partial_{\bar z}g_{k+1}=0\quad\text{for }k\in\Z.
\end{equation}
That is, the invariance under the free motion is a system of $\partial_z$- and $\partial_{\bar z}$-equations for the Fourier coefficients.

As in Section \ref{sec:twistor}, regard $w=\omega v\in\D$ as a vector in $T_x\Omega$ of length less than $1$. Then $g_k(x)e^{ik\theta}\omega^k=g_k(x)w^k$. So the series in Section \ref{sec:twistor} is the power series
\[U(x,w)=\sum_{k\ge0}g_k(x)w^k.\]
Moreover,
\[\left(\omega^2\eta_++\eta_-\right)U=e^{-i\theta}\left(\partial_{\bar z}+w^2\partial_z\right)U.\]
So the equations \eqref{eq:CR} become
\begin{equation}\label{eq:CRflat}
\partial_{\bar w}U=0\quad\text{and}\quad\left(\partial_{\bar z}+w^2\partial_z\right)U=0.
\end{equation}
Indeed, for the power series $U$,
\[\left(\partial_{\bar z}+w^2\partial_z\right)U=\sum_{k\ge0}\left(\partial_{\bar z}g_k+\partial_zg_{k-2}\right)w^k,\quad\text{where }g_{-1}=g_{-2}=0,\]
and these coefficients vanish if and only if \eqref{eq:Fourierflat} holds for $k\ge-1$.

The equations \eqref{eq:CRflat} are the Cauchy--Riemann equations of an explicit complex structure. Let
\[\xi=z-w^2\bar z.\]
For $|w|<1$, the real linear map $\Lambda_w(z)=z-w^2\bar z$ is invertible. So the map $(x,w)\mapsto(w,\xi)$ is a diffeomorphism from $\Omega\times\D$ onto a domain in $\C^2$. A function on $\Omega\times\D$ is called holomorphic if it is a holomorphic function of $(w,\xi)$.

\begin{lemma}[Holomorphic functions on a flat table]\label{lemma:holo}
A locally bounded function $U(x,w)$ on $\Omega\times\D$ is holomorphic if and only if it satisfies \eqref{eq:CRflat} in the sense of distributions.
\end{lemma}

\begin{proof}
Write $U(x,w)=G(w,\xi)$. Since $\xi$ is holomorphic in $w$, we have
\[\partial_{\bar w}U=\partial_{\bar w}G+\partial_{\bar\xi}G\partial_{\bar w}\bar\xi.\]
Since
\[\left(\partial_{\bar z}+w^2\partial_z\right)\xi=0\quad\text{and}\quad\left(\partial_{\bar z}+w^2\partial_z\right)\bar\xi=1-|w|^4,\]
we have
\[\left(\partial_{\bar z}+w^2\partial_z\right)U=\left(1-|w|^4\right)\partial_{\bar\xi}G.\]
So \eqref{eq:CRflat} is equivalent to
\[\partial_{\bar w}G=\partial_{\bar\xi}G=0.\]
By the elliptic regularity of the Cauchy--Riemann equations, a distribution with these properties is a holomorphic function. This proves the lemma.
\end{proof}

In particular, for fixed $w\in\D$, a holomorphic function $U$ on $\Omega\times\D$ is a holomorphic function of the single variable $\xi$ on the planar domain $\Lambda_w(\Omega)$. We use this in Section \ref{sec:difficulty} to take boundary values at $\partial\Omega$.

\subsection{Steps 1--3 in the interior}\label{sec:interior}
Let $\mu$ be as in Theorem \ref{thm:billiard}. Steps 1--3 in Section \ref{sec:surfaceproof} only use the invariance of $\mu$ in the interior of the table. With the formulas in Section \ref{sec:flat}, they can be carried out by hand.

\textit{Step 1: Removal of the odd Fourier modes.} Recall the velocity reversal $R(x,\theta)=(x,\theta+\pi)$. Since $R\circ a_t=a_{-t}\circ R$, the measure $R_*\mu$ is invariant. Let
\[\nu=\frac12\left(\mu+R_*\mu\right).\]
Then $\nu$ is invariant, $\Pi_*\nu=d\Vol_\Omega$, and $R_*\nu=\nu$. Disintegrate
\[d\nu(x,\theta)=d\nu_x(\theta)\,d\Vol_\Omega(x),\]
where $\nu_x$ is a probability measure on $\Sb$ for almost every $x\in\Omega$. For $k\in\Z$, write
\[c_k(x)=\int_{\Sb}e^{-ik\alpha}\,d\nu_x(\alpha).\]
Then the weight-$k$ Fourier mode of $\nu$ is $c_k(x)e^{ik\theta}$. As in Section \ref{sec:surfaceproof}, we have $c_0=1$, $c_k=0$ for odd $k$, $|c_k|\le1$, and $c_{-k}=\ol{c_k}$. By (i) of Lemma \ref{lemma:weak} with $\phi$ supported in $\Omega\times\Sb$, we have $Q\nu=0$ in $\Dc'(\Omega\times\Sb)$. By \eqref{eq:Fourierflat}, this reads
\[\partial_zc_{k-1}+\partial_{\bar z}c_{k+1}=0\quad\text{in }\Dc'(\Omega)\text{ for }k\in\Z.\]
These equations are used in Steps 2 and 3.

\textit{Step 2: The uniform projection.} Since $c_0=1$, the equation of weight $1$ gives
\[\partial_{\bar z}c_2=-\partial_zc_0=0.\]
That is, $c_2$ is holomorphic in $\Omega$. This is the equation $\eta_-\nu_2=-\eta_+\nu_0=0$ in Step 2 of Section \ref{sec:surfaceproof}.

\textit{Step 3: The Carath\'eodory transform.} Let
\[H(x,w)=1+2\sum_{k\ge1}c_{2k}(x)w^{2k}=\int_{\Sb}\frac{e^{i\alpha}+w}{e^{i\alpha}-w}\,d\nu_x(\alpha).\]
That is, $H(x,\cdot)$ is the Carath\'eodory transform of $\nu_x$, as in Step 3 of Section \ref{sec:surfaceproof}. Since $|c_{2k}|\le1$, the series converges locally uniformly in $w$ with values in $L^\infty(\Omega)$. So $\partial_{\bar w}H=0$, and we may apply $\partial_{\bar z}+w^2\partial_z$ term by term. We obtain
\[\left(\partial_{\bar z}+w^2\partial_z\right)H=2\left(\partial_{\bar z}c_2\right)w^2+2\sum_{k\ge2}\left(\partial_{\bar z}c_{2k}+\partial_zc_{2k-2}\right)w^{2k}=0.\]
Here, the first term vanishes by Step 2, and the other terms vanish by Step 1. By Lemma \ref{lemma:holo}, $H$ is holomorphic on $\Omega\times\D$. This is Lemma \ref{lemma:extension} on the flat table. As in Step 3 of Section \ref{sec:surfaceproof}, the triangle inequality gives
\begin{equation}\label{eq:Hbound}
\Rer H(x,w)\ge\frac{1-|w|}{1+|w|}\quad\text{and}\quad|H(x,w)|\le\frac{1+|w|}{1-|w|}.
\end{equation}
Since $H$ is continuous, redefine $\nu_x$ for every $x\in\Omega$ as the probability measure whose Carath\'eodory transform is $H(x,\cdot)$. Then $x\mapsto\nu_x$ is continuous in the weak-* topology.

\subsection{Step 4: the difficulty at the boundary}\label{sec:difficulty}
As in Step 4 of Section \ref{sec:surfaceproof}, let
\[F=\frac{H-1}{H+1}.\]
Then $F$ is holomorphic on $\Omega\times\D$, $|F|<1$, and $F(x,0)=0$. Write
\[F(x,w)=\sum_{k\ge1}b_k(x)w^k.\]
For $0<r<1$, let
\[F_r(x,\theta)=F\left(x,re^{i\theta}\right).\]
As in the proof of Lemma \ref{lemma:bounded}, choose $r_j\to1$ such that $F_{r_j}$ converges weak-* to some $f\in L^\infty(\Omega\times\Sb)$. Then the Fourier modes of $f$ are $b_k(x)e^{ik\theta}$ for $k\ge1$, and the other modes vanish. By \eqref{eq:CRflat}, we have
\[\partial_{\bar z}b_k+\partial_zb_{k-2}=0\quad\text{for }k\ge1,\quad\text{where }b_0=b_{-1}=0.\]
These are the equations in \eqref{eq:Fourierflat} for the Fourier modes of $f$. So $Qf=0$ in $\Dc'(\Omega\times\Sb)$. That is, $f$ is constant along almost every chord of $\Omega$.

On a surface, this is the invariance of $f$, and the ergodicity concludes the proof of Lemma \ref{lemma:bounded}. On a billiard table, by (ii) of Lemma \ref{lemma:weak}, the invariance also requires that $f$ takes the same value on the incoming and the outgoing chords at almost every reflection. This matching condition cannot be derived directly. Indeed, the reflection reverses the orientation of the fibers. For $k\ge1$,
\[e^{ik\tau_x(\theta)}=\left(-e^{2i\psi(x)}\right)^ke^{-ik\theta}.\]
So $\theta\mapsto f(x,\tau_x(\theta))$ only has Fourier modes of negative weight, whereas $f(x,\cdot)$ only has Fourier modes of positive weight. Hence, the matching condition would force $f$ to vanish at the boundary, which is as strong as the conclusion.

We overcome this difficulty in two steps. First, we take the boundary values $F^\partial$ of $F$ on the fibers over $\partial\Omega$. This boundary is not present for surfaces. On the flat table, the boundary values reduce to those of bounded holomorphic functions of one variable; see Lemma \ref{lemma:bv}. Second, the invariance of $\nu$ at $\partial\Omega$ gives a weaker symmetry, which is compatible with the holomorphic structure. Write
\[\tau_x(w)=-e^{2i\psi(x)}\bar w\quad\text{for }x\in\partial\Omega\setminus\Sigma\text{ and }w\in\D\]
for the reflection of $\D$ across the tangent line. Then $\tau_x(e^{i\theta})=e^{i\tau_x(\theta)}$, which is consistent with the reflection of directions. Since $\tau_x$ is anti-holomorphic, the function $w\mapsto\ol{F^\partial(x,\tau_x(w))}$ is holomorphic. In Lemma \ref{lemma:sym}, we show that it equals $F^\partial(x,\cdot)$. Taking the real parts, $\Rer f$ satisfies the matching condition. So $\Rer f$ is invariant; see Proposition \ref{prop:Ref}. Since $f$ only has Fourier modes of positive weight, it is determined by $\Rer f$.

For $x\in\partial\Omega\setminus\Sigma$ and $\ve>0$, write
\[x^\ve=x-\ve e^{i\psi(x)}\]
for the points on the inner normal.

\begin{lemma}[Boundary values]\label{lemma:bv}
For almost every $x\in\partial\Omega$, the limits
\[H^\partial(x,w)=\lim_{\ve\to0^+}H\left(x^\ve,w\right)\quad\text{and}\quad F^\partial(x,w)=\lim_{\ve\to0^+}F\left(x^\ve,w\right)\]
exist locally uniformly in $w\in\D$. They are holomorphic in $w$, and
\[F^\partial=\frac{H^\partial-1}{H^\partial+1}.\]
Moreover, $\nu_{x^\ve}$ converges weak-* to a probability measure $\nu^\partial_x$ whose Carath\'eodory transform is $H^\partial(x,\cdot)$.
\end{lemma}

\begin{proof}
\textit{Step 1: A fixed $w$.} Fix $w\in\D$. As in the proof of Lemma \ref{lemma:holo}, write
\[F(y,w)=G\left(w,\Lambda_w(y)\right)\quad\text{for }y\in\Omega.\]
Then $G(w,\cdot)$ is a holomorphic function on $\Lambda_w(\Omega)$ with $|G|<1$. So the problem reduces to the boundary values of a bounded holomorphic function of one variable.

Let $x\in\partial\Omega\setminus\Sigma$. The points $x^\ve$ lie on the inner normal at $x$. Since $\Lambda_w$ is a real linear isomorphism, the points $\Lambda_w(x^\ve)$ lie on a segment which is transversal to $\Lambda_w(\partial\Omega)$ at $\Lambda_w(x)$. So they approach $\Lambda_w(x)$ non-tangentially. Moreover, $\Lambda_w$ maps the null sets of $\partial\Omega$ to the null sets of $\Lambda_w(\partial\Omega)$. Hence, it suffices to show that $G(w,\cdot)$ has non-tangential limits at almost every point of $\Lambda_w(\partial\Omega\setminus\Sigma)$.

This is a local statement. Near each point of $\Lambda_w(\partial\Omega\setminus\Sigma)$, the curve $\Lambda_w(\partial\Omega)$ is $C^2$. So an open arc around this point lies on the boundary of a Jordan domain $D\subset\Lambda_w(\Omega)$ with $C^2$ boundary. By the Kellogg--Warschawski theorem, a Riemann map $\D\to D$ extends to a $C^1$ diffeomorphism from $\ol\D$ onto $\ol D$; see Pommerenke \cite[Chapter 3]{P}. So it preserves the non-tangential approach and the null sets. Hence, the claim on this arc follows from Fatou's theorem on $\D$; see Garnett \cite[Chapter I]{G}. Countably many such arcs cover $\Lambda_w(\partial\Omega\setminus\Sigma)$.

Therefore, $F(x^\ve,w)$ converges as $\ve\to0^+$ for almost every $x\in\partial\Omega$. By \eqref{eq:Hbound}, $H(\cdot,w)$ is also bounded. So the same holds for $H$.

\textit{Step 2: All $w$.} Fix a countable dense subset of $\D$. By Step 1, there is a null set $N\subset\partial\Omega$ such that $F(x^\ve,w)$ and $H(x^\ve,w)$ converge for $x\in\partial\Omega\setminus N$ and $w$ in this subset.

For $y\in\Omega$, the functions $F(y,\cdot)$ and $H(y,\cdot)$ are holomorphic on $\D$, since the fiber $\{y\}\times\D$ is a complex curve. By \eqref{eq:Hbound} and $|F|<1$, they are locally bounded on $\D$, uniformly in $y$. So Vitali's theorem applies. For $x\in\partial\Omega\setminus N$, the functions $F(x^\ve,\cdot)$ and $H(x^\ve,\cdot)$ converge locally uniformly on $\D$ to holomorphic functions $F^\partial(x,\cdot)$ and $H^\partial(x,\cdot)$. By \eqref{eq:Hbound}, $\Rer H^\partial>0$. So
\[F^\partial=\frac{H^\partial-1}{H^\partial+1}.\]

\textit{Step 3: The conditional measures.} The Taylor coefficients of $H(x^\ve,\cdot)$ are $1$ and $2c_{2k}(x^\ve)$ for $k\ge1$. Up to the factor $2$, they are the Fourier coefficients of $\nu_{x^\ve}$, whose odd Fourier coefficients vanish. By Step 2, these coefficients converge for $x\in\partial\Omega\setminus N$. The probability measures on $\Sb$ are weak-* compact, and they are determined by their Fourier coefficients. So $\nu_{x^\ve}$ converges weak-* to a probability measure $\nu^\partial_x$. Its Carath\'eodory transform is $H^\partial(x,\cdot)$.
\end{proof}

The next lemma is a flux balance at the boundary. The mass of $\nu$ which flows into $\partial\Omega$ equals the mass which flows out after the reflection. Since the reflection is anti-holomorphic, this symmetry of $\nu^\partial_x$ becomes a conjugate symmetry of $F^\partial$.

\begin{lemma}[Reflection symmetry]\label{lemma:sym}
For almost every $x\in\partial\Omega$, the measure $\nu^\partial_x$ is invariant under $\tau_x$. Consequently,
\[F^\partial\left(x,\tau_x(w)\right)=\ol{F^\partial(x,w)}\quad\text{for }w\in\D.\]
\end{lemma}

\begin{proof}
Let $\phi\in C_\tau$. Consider the bounded continuous vector field
\[Y(y)=\int_{\Sb}\phi(y,\theta)(\cos\theta,\sin\theta)\,d\nu_y(\theta)\quad\text{for }y\in\Omega.\]
Since $\phi$ vanishes near $\Sigma\times\Sb$, the vector field $Y$ vanishes near $\Sigma$. For $\phi_0\in C_0^\infty(\Omega)$, the invariance of $\nu$ gives
\[\int Q\left(\phi_0\phi\right)\,d\nu=0.\]
So
\[\Div Y(y)=\int_{\Sb}Q\phi(y,\theta)\,d\nu_y(\theta)\in L^\infty(\Omega).\]
Let
\[\Omega^\ve=\left\{y\in\Omega:\dist(y,\partial\Omega)>\ve\right\}.\]
Fix $\chi\in C^\infty(\R)$ with $\chi=0$ on $(-\infty,1]$ and $\chi=1$ on $[2,\infty)$. For $\ve>0$, let
\[\chi_\ve(y)=\chi\left(\frac{\dist(y,\partial\Omega)}{\ve}\right)\quad\text{for }y\in\Omega.\]
Then $\chi_\ve$ is Lipschitz and supported in $\ol{\Omega^\ve}$. So
\[\int_\Omega\chi_\ve\Div Y\,dy=-\int_\Omega\nabla\chi_\ve\cdot Y\,dy.\]
As $\ve\to0$, the left-hand side converges to
\[\int_\Omega\Div Y\,dy=|\Omega|\int Q\phi\,d\nu\]
by the dominated convergence theorem. On the right-hand side, $\nabla\chi_\ve$ is supported in $\ol{\Omega^\ve}\setminus\Omega^{2\ve}$. Since $Y$ vanishes near $\Sigma$, for small $\ve$, the points of this set in the support of $Y$ have the form $y=x^\delta$ with $x\in\partial\Omega\setminus\Sigma$ and $\ve\le\delta\le2\ve$. In the coordinates $(x,\delta)$, we have
\[\nabla\chi_\ve(y)=-\frac1\ve\chi'\left(\frac\delta\ve\right)e^{i\psi(x)}\quad\text{and}\quad dy=\left(1+O(\delta)\right)\,dsd\delta.\]
By Lemma \ref{lemma:bv}, for almost every $x\in\partial\Omega$,
\[\lim_{\delta\to0^+}e^{i\psi(x)}\cdot Y\left(x^\delta\right)=\int_{\Sb}\phi(x,\theta)\cos\left(\theta-\psi(x)\right)\,d\nu^\partial_x(\theta).\]
Moreover,
\[\int_\ve^{2\ve}\frac1\ve\chi'\left(\frac\delta\ve\right)\,d\delta=1.\]
Since $Y$ is bounded, the dominated convergence theorem gives
\[|\Omega|\int Q\phi\,d\nu=\int_{\partial\Omega}\int_{\Sb}\phi(x,\theta)\cos\left(\theta-\psi(x)\right)\,d\nu^\partial_x(\theta)ds(x).\]
The left-hand side vanishes by (i) of Lemma \ref{lemma:weak}. Take $\phi$ as in the proof of (ii) of Lemma \ref{lemma:weak}. Then the measure
\[\cos\left(\theta-\psi(x)\right)\,d\nu^\partial_x(\theta)\]
is anti-symmetric under $\tau_x$ for almost every $x$. By \eqref{eq:cos}, $\nu^\partial_x$ is symmetric under $\tau_x$ away from the two tangent directions, which are fixed by $\tau_x$. So $\nu^\partial_x$ is invariant under $\tau_x$.

Let $\sigma$ be a real measure on $\Sb$ which is invariant under $\tau_x$. The substitution $\alpha\mapsto\tau_x(\alpha)$ gives
\[\int_{\Sb}\frac{e^{i\alpha}+\tau_x(w)}{e^{i\alpha}-\tau_x(w)}\,d\sigma(\alpha)=\ol{\int_{\Sb}\frac{e^{i\alpha}+w}{e^{i\alpha}-w}\,d\sigma(\alpha)}.\]
So
\[H^\partial\left(x,\tau_x(w)\right)=\ol{H^\partial(x,w)}\quad\text{for }w\in\D.\]
Since $F^\partial=(H^\partial-1)/(H^\partial+1)$, the same holds for $F^\partial$. This proves the lemma.
\end{proof}

\subsection{Step 4: the real part}\label{sec:Ref}
In the proof of Lemma \ref{lemma:bounded}, $F_r$ is tested against $\phi\in C_0^\infty(S\M)$. Here, we test $F_r$ against functions in $C_\tau$. A boundary term at $\partial\Omega$ appears, and Lemma \ref{lemma:sym} cancels its real part.

\begin{prop}[Invariance of the real part]\label{prop:Ref}
The function $\Rer f$ is invariant under the billiard flow.
\end{prop}

\begin{proof}
Let $\phi\in C_\tau$. Fix $0<r<1$. By the second equation in \eqref{eq:CRflat} at $w=re^{i\theta}$, we have
\[\eta_-F_r=-r^2\eta_+F_r.\]
So
\[QF_r=\left(1-r^2\right)\eta_+F_r.\]
By the divergence theorem on $\Omega^\ve\times\Sb$,
\[\int_{\Omega^\ve\times\Sb}F_rQ\phi\,dxd\theta=-\left(1-r^2\right)\int_{\Omega^\ve\times\Sb}\left(\eta_+F_r\right)\phi\,dxd\theta+\int_{\partial\Omega^\ve}\int_{\Sb}F_r\phi\cos\left(\theta-\psi\right)\,d\theta ds.\]
Here, $\ve$ is small, and $\psi$ denotes the angle of the outward unit normal of $\Omega^\ve$. On the support of $\phi$, the boundary $\partial\Omega^\ve$ consists of the points $x^\ve$ with $x\in\partial\Omega\setminus\Sigma$, where this angle is $\psi(x)$. Integrating by parts in $z$ once more, the first term on the right-hand side is bounded by $C(1-r^2)$, where $C$ depends only on $\phi$ and $\Omega$. Letting $\ve\to0$ and using Lemma \ref{lemma:bv}, we obtain
\[\int F_rQ\phi\,dxd\theta=\int_{\partial\Omega}\int_{\Sb}F^\partial\left(x,re^{i\theta}\right)\phi(x,\theta)\cos\left(\theta-\psi(x)\right)\,d\theta ds(x)+O\left(1-r^2\right).\]
Let $r=r_j\to1$. The bounded holomorphic function $F^\partial(x,\cdot)$ has radial limits $F^\partial(x,e^{i\theta})$ for almost every $\theta$. By the dominated convergence theorem,
\[\int fQ\phi\,dxd\theta=\int_{\partial\Omega}\int_{\Sb}F^\partial\left(x,e^{i\theta}\right)\phi(x,\theta)\cos\left(\theta-\psi(x)\right)\,d\theta ds(x).\]
By Lemma \ref{lemma:sym}, $\Rer F^\partial(x,e^{i\theta})$ is symmetric under $\theta\mapsto\tau_x(\theta)$. So is $\phi(x,\theta)$, whereas $\cos(\theta-\psi(x))$ is anti-symmetric by \eqref{eq:cos}. Taking the real parts, we obtain
\[\int\Rer f\cdot Q\phi\,d\Vol_L=0.\]
By (ii) of Lemma \ref{lemma:weak}, $\Rer f$ is invariant. This completes the proof.
\end{proof}

We now complete Step 4. By Proposition \ref{prop:Ref} and the ergodicity of the billiard flow, $\Rer f$ is constant. In Lemma \ref{lemma:bounded}, the ergodicity forces $f$ to be constant. Here, it only forces $\Rer f$ to be constant, but this suffices. Indeed, for $k\ge1$, the weight-$k$ Fourier mode of $\Rer f=(f+\ol f)/2$ is
\[\frac12b_k(x)e^{ik\theta}.\]
Hence, $b_k=0$ for $k\ge1$, and so $F=0$. That is, $H=1$. So $c_{2k}=0$ for $k\ge1$, and $\nu=d\Vol_L$.

\textit{Step 5: Extremality.} This is the same as in Section \ref{sec:surfaceproof}. Notice that
\[d\Vol_L=\frac12\left(\mu+R_*\mu\right),\]
where $\mu$ and $R_*\mu$ are invariant probability measures. Since $d\Vol_L$ is ergodic, it is an extreme point of the convex set of invariant probability measures. So $\mu=R_*\mu=d\Vol_L$. This completes the proof of Theorem \ref{thm:billiard}.

%%%%%%%%%%%%%%%%%%%%%%%%%%%%%%%%%%%%%%%%%
\section{Higher dimensions}\label{sec:higherproof}
In this section, we prove Theorem \ref{thm:entropy}. Throughout the section, let $\M$ be a compact connected Riemannian manifold of dimension $n\ge2$. Suppose that the sectional curvature of $\M$ is at most $-K^2$ for some $K>0$. Then the geodesic flow is Anosov; see Anosov \cite{Ano}. Denote by $d$ the distance on $S\M$ induced by the Sasaki metric; see Sakai \cite[Chapter II, Section 4]{Sak}. Write $\dist_\M$ for the Riemannian distance on $\M$. Since $\Pi$ is a Riemannian submersion for the Sasaki metric, it is $1$-Lipschitz from $(S\M,d)$ to $(\M,\dist_\M)$.

For $z\in S\M$, $\ve>0$, and $m\in\N$, consider the two-sided Bowen ball
\[B_m(z,\ve)=\left\{y\in S\M:d\left(a_ky,a_kz\right)<\ve\text{ for }k\in\Z\text{ with }|k|\le m\right\}.\]
The key estimate is that the projection of $B_m(z,\ve)$ to $\M$ is thin.

\begin{lemma}[Projections of Bowen balls]\label{lemma:Bowen}
There are $\ve_0=\ve_0(\M,K)>0$ and $C=C(\M,K)>0$ such that
\[\Vol_\M\left(\Pi\left(B_m(z,\ve)\right)\right)\le C\ve^ne^{-(n-1)Km}\quad\text{for }z\in S\M,\ 0<\ve\le\ve_0,\text{ and }m\in\N.\]
\end{lemma}

\begin{proof}
\textit{Step 1: Local product coordinates.} Denote by $W^u(z)$ and $W^s(z)$ the strong unstable manifold and strong stable manifold at $z\in S\M$, respectively. Write $d^u$ and $d^s$ for the distances along these manifolds. Since the sectional curvature is at most $-K^2$, the comparison estimates for the stable and unstable Jacobi fields, together with the compactness of $\M$, give $c=c(\M,K)>0$ such that for $t\ge0$,
\[d^u\left(a_ty,a_tz\right)\ge ce^{Kt}d^u(y,z)\quad\text{for }y\in W^u(z),\]
and
\[d^s\left(a_{-t}y,a_{-t}z\right)\ge ce^{Kt}d^s(y,z)\quad\text{for }y\in W^s(z).\]
By the local product structure of the Anosov flow, there are $\ve_0,C_1>0$ with the following property. If $d(y,z)<\ve_0$, then there are $s\in\R$ and $w\in W^u(z)$ with $a_sy\in W^s(w)$ and
\[|s|+d^u(w,z)+d^s\left(a_sy,w\right)\le C_1d(y,z);\]
see Katok--Hasselblatt \cite[Section 18.4]{KH}.

The flow preserves the strong stable and unstable foliations. In particular, $a_{t+s}y\in W^s(a_tw)$ and $a_tw\in W^u(a_tz)$ for $t\in\R$. Let $L\ge1$ be a uniform Lipschitz constant for $a_t$ with $|t|\le1$. For $y\in B_m(z,\ve)$, the integer-time condition gives
\[d(a_ty,a_tz)\le L\ve\quad\text{for }|t|\le m.\]
Decreasing $\ve_0$ if necessary, the local product representation applies and is unique for points at distance at most $L\ve_0$. So by continuity in $t$, the point $a_tw$ and the same shift $s$ give the local product representation of $a_ty$ with respect to $a_tz$ for $|t|\le m$. Hence,
\[d^u(a_tw,a_tz)+d^s(a_{t+s}y,a_tw)\le C_1L\ve\quad\text{for }|t|\le m.\]
Decreasing $\ve_0$ further, we may assume that $C_1\ve_0\le1$.

\textit{Step 2: Projection to $\M$.} Let $0<\ve\le\ve_0$. Suppose that $y\in B_m(z,\ve)$. By Step 1, for $|k|\le m$,
\[d^u\left(a_kw,a_kz\right)\le C_1L\ve\quad\text{and}\quad d^s\left(a_{k+s}y,a_kw\right)\le C_1L\ve.\]
Taking $k=m$ in the first inequality and $k=-m$ in the second, the expansion in Step 1 gives
\[d^u(w,z)\le c^{-1}C_1L\ve e^{-Km}\quad\text{and}\quad d^s\left(a_sy,w\right)\le c^{-1}C_1L\ve e^{-Km}.\]
So
\[d\left(a_sy,z\right)\le C\ve e^{-Km}\quad\text{with }|s|\le C_1\ve\le1.\]
Since $\Pi$ is $1$-Lipschitz and $a_{-s}$ is Lipschitz uniformly for $|s|\le1$, we obtain
\[\dist_\M\left(\Pi(y),\Pi\left(a_{-s}z\right)\right)\le d\left(y,a_{-s}z\right)\le C\ve e^{-Km}.\]
That is, $\Pi(B_m(z,\ve))$ lies in the $C\ve e^{-Km}$-neighborhood of the geodesic segment
\[\left\{\Pi\left(a_sz\right):|s|\le C_1\ve\right\};\]
see Figure \ref{fig:tube}. The volume of this neighborhood is at most $C\ve(\ve e^{-Km})^{n-1}$. This proves the lemma.
\begin{figure}[!ht]
\centering
\begin{tikzpicture}
\draw[thick] (-5.8,0) -- (5.8,0);
\draw[dashed] (-4.3,.7) -- (4.3,.7);
\draw[dashed] (-4.3,-.7) -- (4.3,-.7);
\draw[dashed] (-4.3,.7) arc (90:270:.7);
\draw[dashed] (4.3,-.7) arc (-90:90:.7);
\fill (0,0) circle (.06);
\draw (0,0) -- (.7,1.2) node[above right] {$\Pi(z)$};
\draw[|-|] (-4.3,-1.45) -- (4.3,-1.45);
\node at (0,-1.95) {$O(\ve)$};
\draw[<->] (4.0,0) -- (4.0,.7);
\node[anchor=west] at (4.0,1.45) {$O(\ve e^{-Km})$};
\node at (0,-2.65) {$\Pi\left(B_m(z,\ve)\right)$};
\end{tikzpicture}
\caption{The projection $\Pi\left(B_m(z,\ve)\right)$ lies in a tube of length $O(\ve)$ and radius $O(\ve e^{-Km})$ around the geodesic through $\Pi(z)$.}\label{fig:tube}
\end{figure}
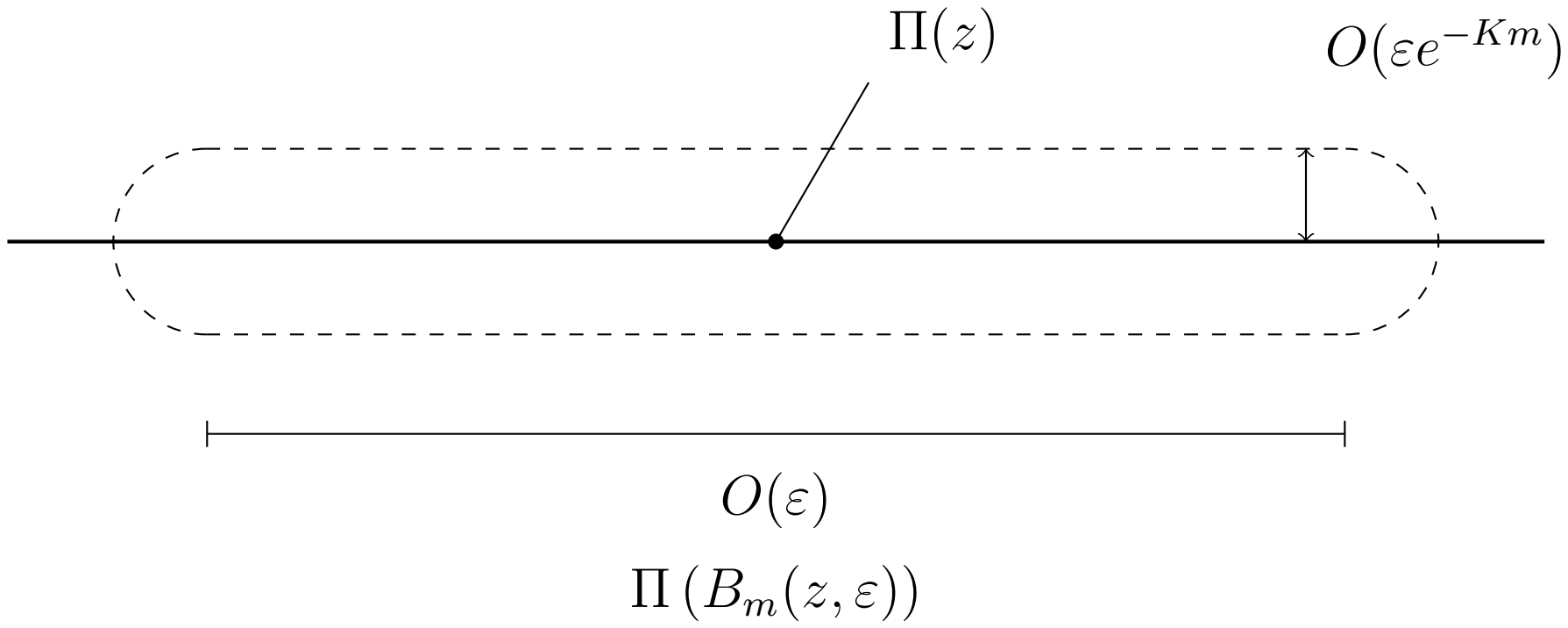
\end{proof}

\begin{proof}[Proof of Theorem \ref{thm:entropy}]
Let $\nu$ be an invariant probability measure with $\nu\le C_2\mu$ for some $C_2>0$. We show that
\[h_\nu(a_1)\ge\frac{(n-1)K}{2}.\]
Fix $0<\ve\le\ve_0$, where $\ve_0$ is from Lemma \ref{lemma:Bowen}. Since $\Pi_*\nu\le C_2\Pi_*\mu=C_2\,d\Vol_\M$, Lemma \ref{lemma:Bowen} gives
\[\nu\left(B_m(z,\ve)\right)\le C_2\Vol_\M\left(\Pi\left(B_m(z,\ve)\right)\right)\le CC_2\ve^ne^{-(n-1)Km}\quad\text{for }z\in S\M\text{ and }m\in\N.\]
Let $\Pc$ be a finite measurable partition of $S\M$ into sets of diameter less than $\ve$. For $m\in\N$, write
\[\Pc^{(2m)}=\bigvee_{k=0}^{2m}a_{-k}\Pc.\]
Let $A\in\Pc^{(2m)}$ and $y\in A$. Then
\[d\left(a_ky',a_ky\right)<\ve\quad\text{for }y'\in A\text{ and }0\le k\le2m.\]
That is,
\[A\subset a_{-m}\left(B_m(a_my,\ve)\right).\]
By the invariance of $\nu$,
\[\nu(A)\le\nu\left(B_m(a_my,\ve)\right)\le CC_2\ve^ne^{-(n-1)Km}.\]
Therefore,
\[-\sum_{A\in\Pc^{(2m)}}\nu(A)\log\nu(A)\ge(n-1)Km-\log\left(CC_2\ve^n\right).\]
Dividing by $2m+1$ and letting $m\to\infty$, we obtain
\[h_\nu(a_1)\ge h_\nu(a_1,\Pc)\ge\frac{(n-1)K}{2}.\]
This proves the claim.

To prove the estimate on each ergodic component, suppose that the invariant set
\[A_0=\left\{z\in S\M:h_{\mu_z}(a_1)<\frac{(n-1)K}{2}\right\}\]
has $\mu(A_0)>0$. Then $\nu=\mu(A_0)^{-1}\mu|_{A_0}$ is an invariant probability measure with $\nu\le\mu(A_0)^{-1}\mu$, and its ergodic components are $\mu_z$ for $\nu$-almost every $z$. By the ergodic decomposition of the entropy,
\[\frac{(n-1)K}{2}\le h_\nu(a_1)=\int_{S\M}h_{\mu_z}(a_1)\,d\nu(z)<\frac{(n-1)K}{2},\]
a contradiction. This completes the proof.
\end{proof}

\appendix
\section{Hyperbolic toral automorphisms}\label{sec:toral}
In this appendix, we show that the analogue of Problem \ref{prob:uniform} has a negative answer for every hyperbolic automorphism of the two-dimensional torus $\T^2=\R^2/\Z^2$. Let $dz$ be the Lebesgue measure on $\T^2$ with coordinates $z=(x,y)$. Write $\Pi(x,y)=x$ for the projection to the position. Denote by $dx$ the Lebesgue measure on $\T$. Let $M\in\GL(2,\Z)$ be hyperbolic. That is, $M$ has no eigenvalues of modulus $1$. It induces an automorphism of $\T^2$, also denoted by $M$.

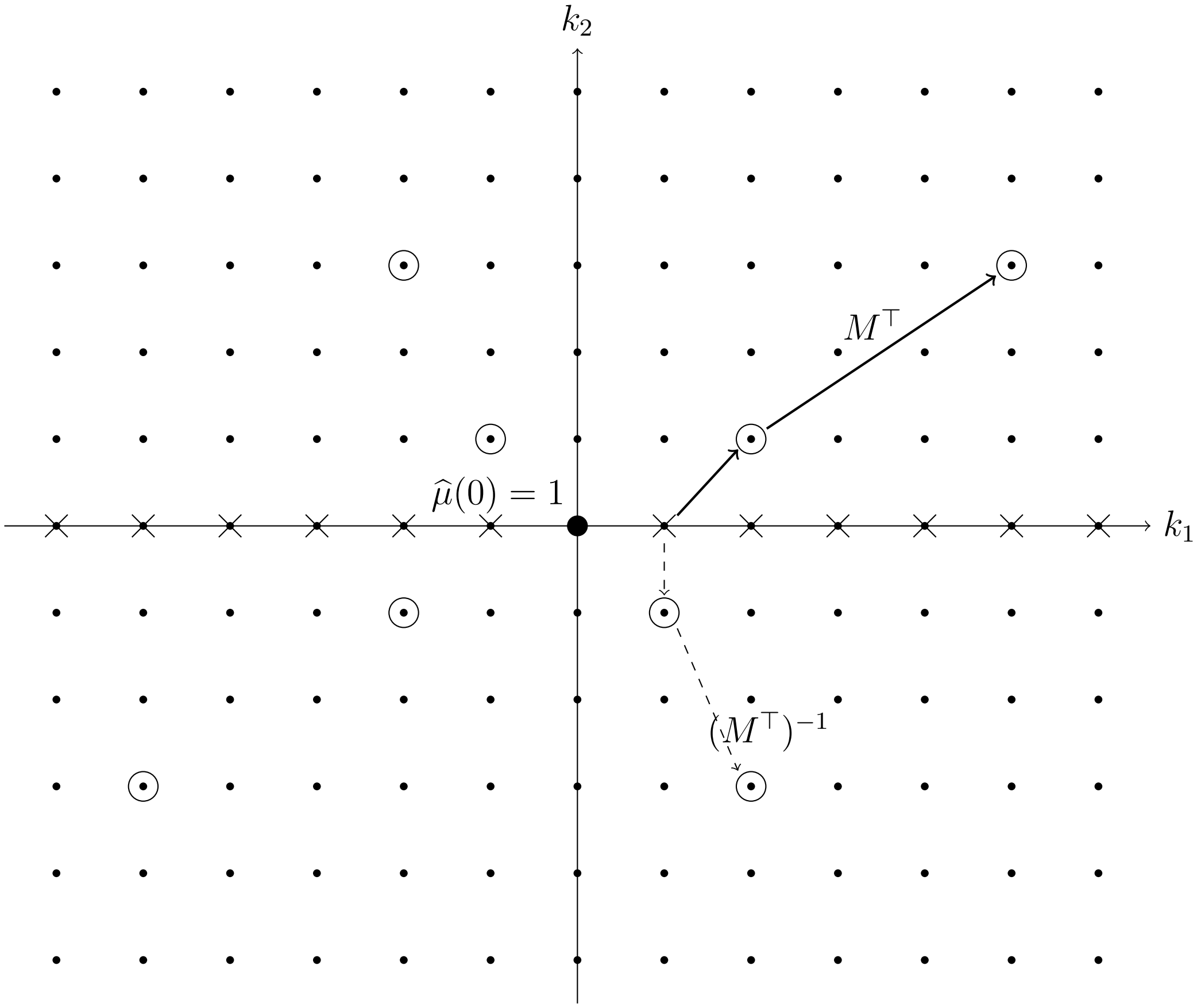
\begin{figure}[!ht]
\centering
\begin{tikzpicture}
\foreach \j in {-6,...,6} {
\foreach \k in {-5,...,5} \fill (\j,\k) circle (.045);
}
\draw[->] (-6.6,0) -- (6.6,0) node[right] {$k_1$};
\draw[->] (0,-5.5) -- (0,5.5) node[above] {$k_2$};
\foreach \j in {-6,-5,-4,-3,-2,-1,1,2,3,4,5,6} {
\draw (\j-.13,-.13) -- (\j+.13,.13);
\draw (\j-.13,.13) -- (\j+.13,-.13);
}
\foreach \j/\k in {2/1,5/3,1/-1,2/-3,-2/-1,-5/-3,-1/1,-2/3} \draw (\j,\k) circle (.17);
\draw[->,thick] (1.15,.12) -- (1.85,.88);
\draw[->,thick] (2.18,1.12) -- (4.82,2.88);
\draw[->,dashed] (1,-.2) -- (1,-.8);
\draw[->,dashed] (1.15,-1.18) -- (1.85,-2.82);
\fill (0,0) circle (.12);
\node[above left] at (0,0) {$\widehat\mu(0)=1$};
\node[above] at (3.4,2) {$M^\top$};
\node[below] at (2.2,-2) {$(M^\top)^{-1}$};
\end{tikzpicture}
\caption{A finite window of the Fourier lattice $\Z^2$ for \eqref{eq:cat}. Uniform projection of $\mu$ forces $\widehat\mu(le_1)=0$ for every $l\ne0$ (crosses). Invariance forces $\widehat\mu((M^\top)^mle_1)=0$ for all $m\in\Z$ (selected additional frequencies are circled). Solid and dashed arrows indicate selected forward and backward iterates of $e_1$, respectively.}\label{fig:cat-frequencies}
\end{figure}

\begin{ex}
Let
\begin{equation}\label{eq:cat}
M=\begin{pmatrix}
2 & 1\\
1 & 1
\end{pmatrix}.
\end{equation}
Figures \ref{fig:cat-frequencies} and \ref{fig:cat-directions} illustrate the construction for this matrix.
\end{ex}

We construct a probability measure $\mu$ on $\T^2$ such that
\begin{equation}\label{eq:uniformcat}
M_*\mu=\mu,\quad\Pi_*\mu=dx,\quad\text{and}\quad\mu\perp dz.
\end{equation}
For a finite measure $\rho$ on $\T^2$, write
\[\widehat\rho(k)=\int_{\T^2}e^{-2\pi i\langle k,z\rangle}\,d\rho(z)\quad\text{for }k\in\Z^2.\]
We can formulate the construction in the frequency space. Let $e_1=(1,0)$. Notice that
\[\widehat{M_*\rho}(k)=\widehat\rho\left(M^\top k\right)\quad\text{and}\quad\widehat{\Pi_*\rho}(l)=\widehat\rho\left(le_1\right).\]
Hence, the first two conditions on $\mu$ in \eqref{eq:uniformcat} are equivalent to
\[\widehat\mu\left(M^\top k\right)=\widehat\mu(k)\quad\text{for }k\in\Z^2\quad\text{and}\quad\widehat\mu\left(le_1\right)=0\quad\text{for }l\in\Z\setminus\{0\}.\]
In particular, $\widehat\mu$ vanishes at the frequencies
\begin{equation}\label{eq:Me1}
\left(M^\top\right)^mle_1\quad\text{for }m\in\Z\text{ and }l\in\Z\setminus\{0\}.
\end{equation}
Figure \ref{fig:cat-frequencies} illustrates them for \eqref{eq:cat}. Thus, we seek a probability measure which is invariant under $M$, whose Fourier coefficients vanish at the frequencies in \eqref{eq:Me1}, and which is singular with respect to $dz$.

\begin{rmk}
The construction also gives invariant probability measures with uniform projection which differ from the Lebesgue measure and have entropy arbitrarily close to the maximal entropy. Indeed, for $0<\theta<1$, let
\[\mu_\theta=(1-\theta)\,dz+\theta\mu.\]
Then $\mu_\theta$ is invariant under $M$, satisfies $\Pi_*\mu_\theta=dx$, and differs from $dz$. Since the entropy is affine,
\[h_{\mu_\theta}(M)=(1-\theta)h_{dz}(M)+\theta h_\mu(M)\ge(1-\theta)h_{dz}(M),\]
where $h_{dz}(M)$ is the maximal entropy of $M$. Thus, $h_{\mu_\theta}(M)\to h_{dz}(M)$ as $\theta\to0$.
\end{rmk}

\begin{rmk}[Outline of the construction]
We use a Riesz product; see Zygmund \cite[Section V.7]{Zy}. Given a lacunary sequence of frequencies $w_j$, the Riesz product is a singular probability measure. Its Fourier coefficients are supported on the sums
\[\sum_j\delta_jw_j\quad\text{with }\delta_j\in\{-1,0,1\},\]
and they are given explicitly. We take the frequencies $w_j$ along one orbit of $(M^\top)^p$ for some $p\in\N$. Then $(M^\top)^p$ shifts the frequencies, so the Riesz product is invariant under $M^p$. Its average over $M^m$ with $0\le m<p$ is invariant under $M$.

It remains to keep these sums away from the frequencies in \eqref{eq:Me1}. We use the coordinates given by the eigenvectors of $M^\top$, in which $M^\top$ stretches one coordinate and shrinks the other. In these coordinates, the frequencies in \eqref{eq:Me1} lie on lines through the origin whose slopes form a discrete set. On the other hand, such a sum is dominated by its term of the largest index in the expanding coordinate, and by its term of the smallest index in the contracting coordinate. This determines its slope up to a small factor. So a suitable choice of $w_0$ and $p$ keeps these slopes away from the slopes of the frequencies in \eqref{eq:Me1}.
\end{rmk}

\subsection{Eigencoordinates}
\begin{figure}[!ht]
\centering
\begin{tikzpicture}
\draw[->] (-6.6,0) -- (6.6,0) node[right] {$k_1$};
\draw[->] (0,-5.5) -- (0,5.5) node[above] {$k_2$};
\draw[dashed] (-6.4,0) -- (6.4,0);
\draw[red,dashed] (-6.4,-3.2) -- (6.4,3.2);
\draw[red,dashed] (-6.4,-3.84) -- (6.4,3.84);
\draw[green!60!black,dashed] (-5.3,5.3) -- (5.3,-5.3);
\draw[green!60!black,dashed] ({-5.3*2/3},5.3) -- ({5.3*2/3},-5.3);
\draw[red,thick] (-6.4,{-6.4*(sqrt(5)-1)/2}) -- (6.4,{6.4*(sqrt(5)-1)/2}) node[above] {$\alpha_-=0$};
\draw[green!60!black,thick] ({-5.3*2/(1+sqrt(5))},5.3) -- ({5.3*2/(1+sqrt(5))},-5.3) node[right] {$\alpha_+=0$};
\node[below] at (-5.6,0) {$m=0$};
\node[red,right] at (6.4,3.2) {$m=1$};
\node[red,left] at (-6.4,-3.84) {$m=2$};
\node[green!60!black,right] at (5.3,-5.3) {$m=-1$};
\node[green!60!black,left] at ({-5.3*2/3},5.3) {$m=-2$};
\fill (0,0) circle (.08);
\end{tikzpicture}
\caption{The eigenlines of $M^\top$ for \eqref{eq:cat} in the frequency space $\R^2$. The expanding eigenline $\{\alpha_-=0\}$ has slope $(\sqrt5-1)/2$, and the contracting eigenline $\{\alpha_+=0\}$ has slope $-(\sqrt5+1)/2$. The dashed lines are $\R(M^\top)^me_1$ for $|m|\le2$. They contain the frequencies in \eqref{eq:Me1}, which the construction avoids. As $m\to\infty$ and $m\to-\infty$, these lines approach the expanding and the contracting eigenlines, respectively.}\label{fig:cat-directions}
\end{figure}
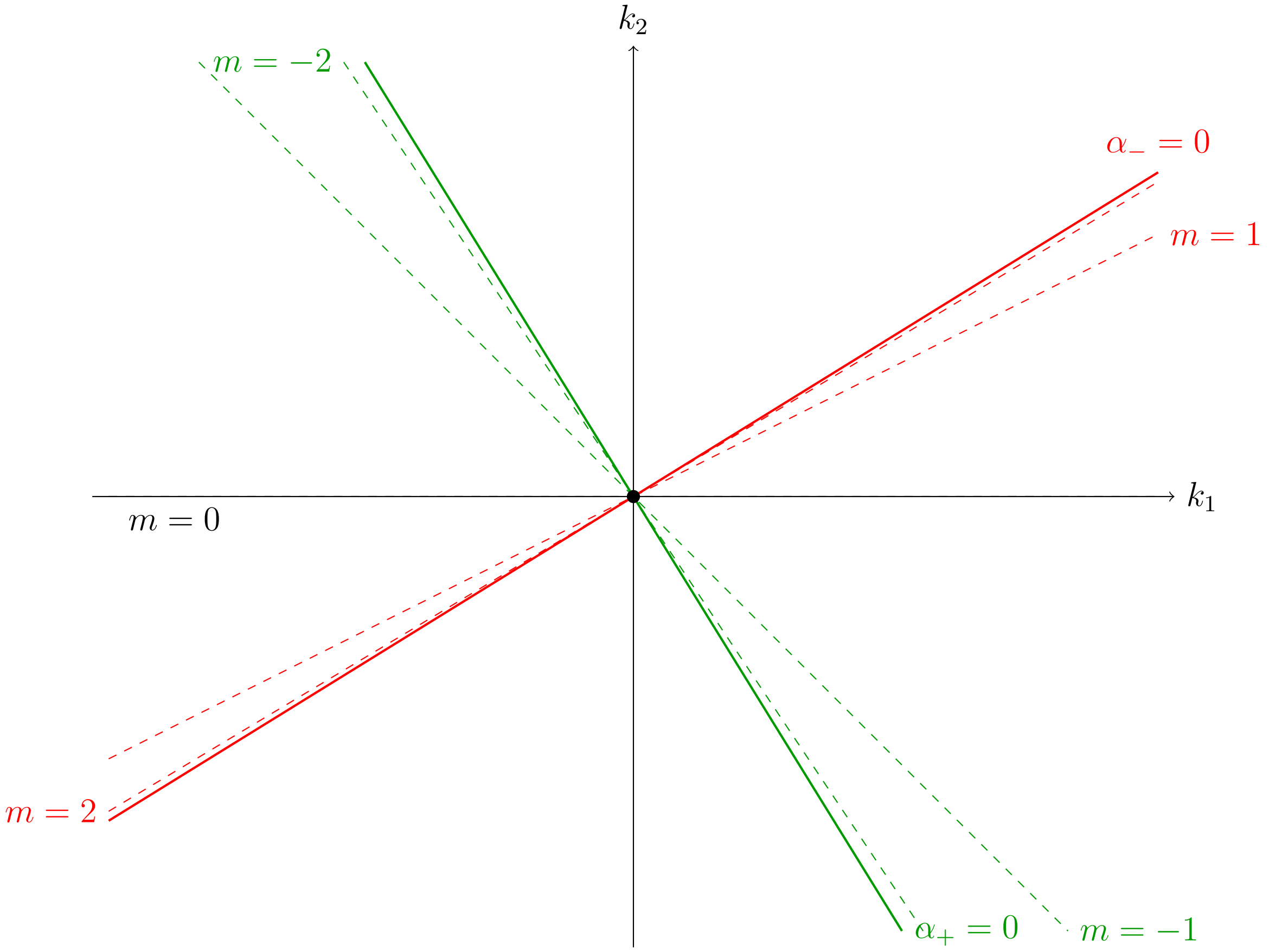

Since $M$ is hyperbolic and $|\det M|=1$, the eigenvalues of $M^\top$ have absolute values $e^\chi$ and $e^{-\chi}$ for some $\chi>0$. If an eigenvalue of $M^\top$ were rational, then it would be an integer dividing $\det M=\pm1$, which is impossible. So the eigenvalues are irrational, and no non-zero vector in $\Z^2$ is an eigenvector of $M^\top$. Denote by $\alpha_+(w)$ and $\alpha_-(w)$ the coordinates of $w\in\R^2$ with respect to eigenvectors of $M^\top$ for the expanding and the contracting eigenvalues, respectively. Then
\[\left|\alpha_\pm\left(M^\top w\right)\right|=e^{\pm\chi}\left|\alpha_\pm(w)\right|\quad\text{for }w\in\R^2.\]
That is, $M^\top$ stretches the coordinate $\alpha_+$ and shrinks the coordinate $\alpha_-$; see Figure \ref{fig:cat-directions}. Moreover, $\alpha_\pm(k)\ne0$ for $k\in\Z^2\setminus\{0\}$.

To locate the lines through the origin, write
\[\Phi(w)=\left|\frac{\alpha_+(w)}{\alpha_-(w)}\right|\quad\text{for }w\in\R^2\text{ with }\alpha_+(w)\alpha_-(w)\ne0.\]
That is, $\Phi(w)$ is the absolute value of the slope of the line $\R w$ in the eigencoordinates. Then $\Phi(tw)=\Phi(w)$ for $t\ne0$, and
\[\Phi\left(M^\top w\right)=e^{2\chi}\Phi(w).\]
Hence, at the frequencies in \eqref{eq:Me1}, the function $\log\Phi$ takes values in the discrete set $\log\Phi(e_1)+2\chi\Z$.

Choose the first frequency $w_0$ such that $\log\Phi(w_0)$ stays away from this discrete set, even after adding multiples of $\chi$.

\begin{lemma}[Initial frequency]\label{lemma:w0}
There is $w_0\in\Z^2\setminus\{0\}$ such that
\[d_0=\dist\left(\log\Phi\left(w_0\right)-\log\Phi\left(e_1\right),\chi\Z\right)>0.\]
\end{lemma}

\begin{proof}
The value $\Phi(w)$ only depends on the line $\R w$. It tends to $\infty$ as $\R w$ approaches the expanding eigenline, and to $0$ as $\R w$ approaches the contracting eigenline. Since the lines through the points of $\Z^2$ are dense, $\log\Phi(\Z^2\setminus\{0\})$ is dense in $\R$. So it is not contained in the discrete set $\log\Phi(e_1)+\chi\Z$. This proves the lemma.
\end{proof}

\subsection{Lacunary frequencies}
Fix $w_0$ and $d_0$ as in Lemma \ref{lemma:w0}. Choose $p\in\N$ such that $\beta=e^{p\chi}$ satisfies
\[\beta>3\quad\text{and}\quad\log\frac{\beta}{\beta-2}<d_0.\]
This is possible, since $d_0>0$, and the left-hand side tends to $0$ as $p\to\infty$. Let
\[w_j=\left(M^\top\right)^{pj}w_0\quad\text{for }j\in\Z.\]
Then
\[\left|\alpha_+\left(w_j\right)\right|=\beta^j\left|\alpha_+\left(w_0\right)\right|\quad\text{and}\quad\left|\alpha_-\left(w_j\right)\right|=\beta^{-j}\left|\alpha_-\left(w_0\right)\right|.\]
That is, the frequencies $w_j$ grow geometrically with ratio $\beta$ in the coordinate $\alpha_+$ as $j\to\infty$, and in the coordinate $\alpha_-$ as $j\to-\infty$. We say that $k\in\Z^2$ has a digit expansion if
\[k=\sum_{j\in\Z}\delta_jw_j\quad\text{with }\delta_j\in\{-1,0,1\}\]
and finitely many $\delta_j\ne0$.

\begin{lemma}[Digit expansions]\label{lemma:digits}\leavevmode
\begin{enumerate}[(i)]
\item Let $c_j\in\Z$ with $|c_j|\le2$ for $j\in\Z$. Suppose that at least one and at most finitely many $c_j$ are non-zero. Then
\[\sum_jc_jw_j\ne0.\]
In particular, digit expansions are unique.
\item If $k\ne0$ has a digit expansion, then $k$ is not among the frequencies in \eqref{eq:Me1}.
\end{enumerate}
\end{lemma}

The proof uses that, since $\beta>3$, a sum as in (i) is dominated by its term of the largest index in the coordinate $\alpha_+$, and by its term of the smallest index in the coordinate $\alpha_-$.

\begin{proof}\leavevmode
\begin{enumerate}[(i)]
\item Denote by $J$ the largest $j$ with $c_j\ne0$. Since
\[\sum_{j<J}\beta^j=\frac{\beta^J}{\beta-1},\]
we have
\[\left|\alpha_+\left(\sum_jc_jw_j\right)\right|\ge\left|\alpha_+\left(w_0\right)\right|\left(\beta^J-2\sum_{j<J}\beta^j\right)=\left|\alpha_+\left(w_0\right)\right|\beta^J\cdot\frac{\beta-3}{\beta-1}>0.\]
So the sum does not vanish. The difference of two digit expansions of the same frequency has coefficients $c_j\in\Z$ with $|c_j|\le2$, and it vanishes. So all its coefficients vanish, and the two expansions coincide. This proves (i).

\item Let
\[k=\sum_j\delta_jw_j\]
be the digit expansion. Denote by $J_+$ and $J_-$ the largest and the smallest $j$ with $\delta_j\ne0$, respectively. As in (i),
\[\left|\alpha_+(k)\right|=r_+\beta^{J_+}\left|\alpha_+\left(w_0\right)\right|\quad\text{and}\quad\left|\alpha_-(k)\right|=r_-\beta^{-J_-}\left|\alpha_-\left(w_0\right)\right|,\quad\text{where }\frac{\beta-2}{\beta-1}\le r_\pm\le\frac{\beta}{\beta-1}.\]
Hence,
\[\Phi(k)=\frac{r_+}{r_-}\beta^{J_++J_-}\Phi\left(w_0\right),\quad\text{where }\frac{\beta-2}{\beta}\le\frac{r_+}{r_-}\le\frac{\beta}{\beta-2}.\]
Suppose that $k$ is among the frequencies in \eqref{eq:Me1}. Then $\log\Phi(k)\in\log\Phi(e_1)+2\chi\Z$. Since $\log\beta=p\chi$, the distance from $\log\Phi(w_0)-\log\Phi(e_1)$ to $\chi\Z$ is at most $\log(\beta/(\beta-2))$. This contradicts the choice of $p$, and proves (ii).
\end{enumerate}
\end{proof}

\subsection{The Riesz product}
Fix $q\in(0,\frac12)$. For $N\in\N$, define
\[R_N(z)=\prod_{j=-N}^N\left(1+qe^{2\pi i\langle w_j,z\rangle}+qe^{-2\pi i\langle w_j,z\rangle}\right)\quad\text{for }z\in\T^2.\]
Since $0<q<\frac12$, each factor is positive. So $R_N$ is positive. Expanding the product, we obtain
\[R_N(z)=\sum_{(\delta_j)_{|j|\le N}\in\{-1,0,1\}^{2N+1}}q^{\sum_j|\delta_j|}e^{2\pi i\left\langle\sum_j\delta_jw_j,z\right\rangle}.\]
By (i) of Lemma \ref{lemma:digits}, distinct choices of $(\delta_j)$ give distinct frequencies, and the frequency vanishes only if all $\delta_j$ vanish. So
\[\int_{\T^2}R_N\,dz=1,\]
that is, $R_N\,dz$ is a probability measure. Moreover, its Fourier coefficient at $k\in\Z^2$ is $q$ to the power of the number of non-zero digits if $k$ has a digit expansion with $\delta_j=0$ for $|j|>N$, and $0$ otherwise. By the uniqueness of digit expansions, this coefficient does not depend on $N$ for large $N$. So $R_N\,dz$ converges weakly to a probability measure $\rho$ on $\T^2$ with
\[\widehat\rho(k)=\begin{cases}q^{\sum_j|\delta_j|}&\text{if }k=\sum_j\delta_jw_j\text{ is a digit expansion},\\0&\text{if }k\text{ has no digit expansion}.\end{cases}\]
Since $(M^\top)^pw_j=w_{j+1}$, the frequency $k$ has a digit expansion if and only if $(M^\top)^pk$ has one, with the same number of non-zero digits. So
\[\widehat\rho\left(\left(M^\top\right)^pk\right)=\widehat\rho(k)\quad\text{for }k\in\Z^2,\]
that is, $(M^p)_*\rho=\rho$. By (ii) of Lemma \ref{lemma:digits},
\[\widehat{\left(M^m\right)_*\rho}(le_1)=\widehat\rho\left(\left(M^\top\right)^mle_1\right)=0\quad\text{for }l\in\Z\setminus\{0\}\text{ and }m\in\Z.\]
That is, $\Pi_*(M^m)_*\rho=dx$ for every $m\in\Z$.

Now construct the probability measure in this appendix. Let
\[\mu=\frac1p\sum_{m=0}^{p-1}\left(M^m\right)_*\rho.\]
Since $(M^p)_*\rho=\rho$, we have $M_*\mu=\mu$. Since $\Pi_*(M^m)_*\rho=dx$ for each $m$, we have $\Pi_*\mu=dx$.

It remains to prove that $\mu\perp dz$. The singularity comes from the lacunarity. Under $\rho$, the functions $\cos(2\pi\langle w_j,z\rangle)$ are uncorrelated and have mean $q$. Under $dz$, they have mean $0$. A law of large numbers separates the two measures. For $N\ge1$, let
\[S_N(z)=\frac1N\sum_{j=0}^{N-1}\cos\left(2\pi\langle w_j,z\rangle\right).\]
By (i) of Lemma \ref{lemma:digits},
\[\widehat\rho\left(\pm w_j\right)=q\quad\text{and}\quad\widehat\rho\left(\pm w_j\pm w_{j'}\right)=q^2\text{ for }j\ne j'.\]
Since $|\widehat\rho(2w_j)|\le1$, we obtain
\[\int_{\T^2}S_N\,d\rho=q\quad\text{and}\quad\int_{\T^2}\left|S_N-q\right|^2\,d\rho\le\frac{1-q^2}{N}\le\frac1N.\]
On the other hand, $w_j\ne\pm w_{j'}$ for $j\ne j'$ by (i) of Lemma \ref{lemma:digits}. So the orthogonality of characters gives
\[\int_{\T^2}S_N\,dz=0\quad\text{and}\quad\int_{\T^2}\left|S_N\right|^2\,dz=\frac1{2N}.\]
By the Borel--Cantelli lemma, $S_{2^l}\to q$ $\rho$-almost everywhere and $S_{2^l}\to0$ Lebesgue-almost everywhere as $l\to\infty$. Hence, $\rho\perp dz$. Since $M$ preserves $dz$, also $(M^m)_*\rho\perp dz$ for each $m$, and so $\mu\perp dz$.

%%%%%%%%%%%%%%%%%%%%%%%%%%%%%%%%%%%%%%%%%
\section*{Acknowledgements}
The author has benefited greatly from discussions with Manfred Einsiedler, Giovanni Forni, Michael Hochman, Fran\c{c}ois Ledrappier, Anthony Quas, Tuomas Sahlsten, Pablo Shmerkin, and Zhiyuan Zhang. He is grateful to them for their valuable insights. He also thanks Ze\'ev Rudnick for comments on an earlier version of this paper, and Giovanni Forni for informing him of the paper \cite{RRW}.

%%%%%%%%%%%%%%%%%%%%%%%%%%%%%%%%%%%%%%%%%
\section*{Use of artificial intelligence}
The author used Claude (Anthropic) to assist with mathematical arguments, exposition, and reference searches. The author takes full responsibility for the content of this paper.

%%%%%%%%%%%%%%%%%%%%%%%%%%%%%%%%%%%%%%%%%

\end{document}